\documentclass[twoside]{siamltex}
\usepackage{latexsym,epsf}
\def\R{\hbox{\rm I\kern-.2em\hbox{\rm R}}}

\newcommand{\norm}[1]{\left\Vert #1 \right\Vert}

\newtheorem{thm}{Theorem}[section]
\newtheorem{lem}[thm]{Lemma} 

\newtheorem{exm}[thm]{\it Example\/\rm } 
\newtheorem{rem}{\it Remark\/\rm }[section]
\def\bmat#1{\left[\matrix{#1}\right]}

\def\AMSMOS{\@abssec{AMS(MOS) subject classifications}}
\def\@abssec#1{\vspace{.05in}\footnotesize \parindent .2in
{\bf #1. }\ignorespaces}

\begin{document}
\title{ORTHONORMAL INTEGRATORS BASED ON HOUSEHOLDER
AND GIVENS TRANSFORMATIONS\thanks{This work 
was supported in part under NSF Grants \#DMS-9973266
and \#DMS-9973393.} }
\author{LUCA DIECI\thanks{School of Mathematics and CDSNS, 
Georgia Institute of Technology, Atlanta, GA 30332-0160 U.S.A..}
\and 
ERIK S. VAN VLECK\thanks{
Dept. of Mathematical \& Computer Sciences, Colorado School
of Mines, Golden, CO 80401 U.S.A..}}

\null\hfill Draft version of \today

\maketitle
\begin{abstract}
We carry further our work \cite{DiVa:ONFact} on 
orthonormal integrators based on Householder and Givens
transformations.   We propose new algorithms and pay
particular attention to appropriate implementation
of these techniques.  We also present a suite of 
{\tt Fortran} codes and provide 
numerical testing to show the efficiency and accuracy of 
our techniques.
\end{abstract}

\begin{AMSMOS} Primary 65L \end{AMSMOS}


\section{Introduction}\label{s1}
In recent times, there has been a lot of interest in techniques
for solving differential equations whose exact solution is
an orthonormal matrix, in such a way that the computed solution
is also orthonormal; \cite{CaIsZa:NumSolIsosp}, \cite{DRV}, 
\cite{DiVa:CompLyaExp}, \cite{DiVa:ONFact}, 
\cite{Diele:NumMetUniDiffSys}, \cite{GoSuOr:StabLyaStabDynSys},
\cite{Higham:TimeStepPresOrth}, \cite{MeRa:NumMethSVD},
are a representative sample of references.
In this paper, we will restrict attention to a particular
class of differential equations with orthonormal solution, the 
{\sl $QR$ equations}. 
These arise in many seemingly unrelated applications
(e.g., see \cite{BGGS}, \cite{Chu:ContRealIterProc}, 
\cite{DeNaTo:ODESymEigenProb}, \cite{DiVa:ONFact},
\cite{Me}), but can all be tracked back to the following setup.

\noindent 1.
We are given the function $A:\ t\to \R^{n\times n}$, 
$A\in C^{k-1},\ k\ge 1$, and consider the associated initial value problem 
for $X\in \R^{n\times p}$, $p\le n$:
$$\dot X\ = \ A(t) X\,,\,\, t\in [t_0,t_f]\ ,\,\,\, X(t_0)=X_0 \,\,\, 
\rm{full}\,\,\, \rm{rank}\,.$$
\noindent 2.  We need the $QR$ factorization of 
the matrix $X$: $X=QR$, $Q, R\in C^{k}$, $Q\in \R^{n\times p},
\ R\in \R^{p\times p}$.  For stability reasons, we
want to find the factors $Q$ and $R$ without first finding $X$.  Once 
$Q$ is known, the equation for $R$ becomes:
\begin{equation}\label{Rflow}
\dot R\ = \ \tilde A R\ ,\,\, R(t_0)=R_0\,,\,\,\, 
\tilde A(t)= (Q^TAQ-Q^T\dot Q) R\,,
\end{equation}
where --since $X=QR$-- $\tilde A$ of (\ref{Rflow}) is upper 
triangular.  Computation of $R$ reduces to a backward 
substitution algorithm for the entries of $R$, coupled
with quadratures,  
which is {\sl conceptually} simple.  Therefore,
in this work we will focus on finding $Q$.  \hfill\break
\noindent 3.  Let $Q_0R_0=X_0$ be a $QR$ factorization of $X_0$.
Differentiate the relation $X=QR$, make use of 
triangularity of $R$, let $S=Q^T\dot Q$, observe that 
$S\in \R^{p\times p}$ must be skew-symmetric, and obtain the
following differential equation for $Q$:
\begin{equation}\label{QRflow}
\dot Q= AQ-QQ^TAQ+QS\,,\,\ 
S_{ij} =\cases{(Q^TAQ)_{ij},& $i> j$ \cr 
                 0,  &  $i=j$ \cr
   -(Q^T AQ)_{ji},& $i<j$\cr }\,\,\,, \ Q(t_0)=Q_0\,.
\end{equation}

\begin{rem}\rm
Different rewritings of (\ref{QRflow}) are possible, in the case
of $p<n$, which present distinct advantages with respect to the given 
form (\ref{QRflow}).  In particular, the following rewriting proposed
by Bridges and Reich in \cite{BridReich:Stiefel} is interesting:
\begin{equation}\label{QR-BR}
\dot Q=[(A-QQ^TA)-(A^T-A^TQQ^T)+QSQ^T]Q\,.
\end{equation}
\end{rem}

\begin{rem}\rm
Triangularity of $\tilde A$ is the key point of a successful and widely 
used technique for computing Lyapunov exponents of the $A$-system
in decreasing order, i.e., the
first $p$ dominant ones; e.g., see \cite{BGGS} and 
\cite{DiRuVa:CompLyaExpContDynSys}.  It is worth pointing out 
that one may not know a priori how many Lyapunov
exponents are needed, i.e., the value of $p$; for example, one may
need to compute all the positive ones\footnote{e.g., this is needed in
the study of ergodic dynamical systems to approximate the so-called
Lyapunov dimension and physical measure, see
\cite{EcRu:ErgTheChao}}.  In this case, it would
be desirable to have techniques which are capable 
(at least, in principle) of increasing the number of exponents to
be calculated without having to restart the entire computation from
scratch, and that can profitably use the computations done so far.
\end{rem}

In the next section, we discuss the design choices we 
face when devising methods to solve (\ref{QRflow}).  Then,
we revisit the setup we proposed in \cite{DiVa:ONFact} where 
we integrate (\ref{QRflow}) by seeking $Q$ as product of elementary
Householder or Givens transformations.  In Section \ref{s3} we discuss
numerical issues which must be confronted when implementing techniques
based on these elementary transformations.  We also put forward new 
formulations in the Householder case.  In Section \ref{s4}, we present 
{\tt FORTRAN} codes we have written for the task, and
we illustrate their performance on a number of 
examples.  Conclusions are in Section \ref{s5}.

We must point out right away that in this work we are exclusively
concerned with the case of coefficient matrices $A(\cdot)$ which
are {\bf dense}; that is, they have no particularly exploitable
structure.  In case $A$ is structured (symmetric, banded, etc.)
some computational savings ought to be possible, and we will
consider these cases in future work.
 
\section{Background and Householder and Givens representations}\label{s2}
Integrating (\ref{QRflow}) is an example of integration on
a manifold.  Quite clearly, we can view the solution $Q$ as a curve 
on the compact
manifold of orthonormal (orthogonal if $p=n$) matrices.  
The dimension of this manifold is $p\frac{2n-p-1}{2}$, which is
therefore the number of degrees of freedom one has to resolve. In general,
a smooth manifold may be parametrized in many different ways, and the
choice of parametrization may turn out to be of utmost importance from 
the numerical point of view.  Clearly, a minimal parametrization requires
$p\frac{2n-p-1}{2}$ parameters.  We notice that a typical construction
in differential geometry is to parametrize a smooth manifold by 
overlapping local charts; this may
be a convenient way to seek the solution of a differential 
equation evolving on a manifold even when the manifold itself
may happen to be globally parametrizable.  Indeed, 
from the numerical point of view, it is important that the parametrization
lead to a stable numerical procedure.  For example,
it is not enough to know that the solution of a differential
equation on a manifold has a global
parametrization to make adopting this parametrization
numerically legitimate.  In particular, if we choose a
representation for the solution $Q$ of (\ref{QRflow}), even 
when it happens to be globally defined, it
may not give us a sound numerical procedure.

Since the solution of (\ref{QRflow})
exists for all times, one may think that direct integration of
the IVP for $Q$ is the right way to proceed.  However, a standard
numerical method used on (\ref{QRflow}) will not deliver a 
solution which is orthonormal at the meshpoints.  
For this reason, several techniques
have been proposed to obtain an orthonormal numerical solution.

\begin{itemize}
\item[[$p=n$]]  In this so-called ``square case'', there are many
competing approaches.  In our opinion, the ones below are the
most natural and/or interesting.  All of these 
choices have been implemented to varying degree of sophistication,
and all present intrinsic implementation difficulties: without entering
in the specifics of these methods, it is enough to state that
none of the methods below is trivial to implement, with the exception
of (ii), and that all them can be implemented so to
have an expense proportional to $O(n^3)$ per step, which is the
expense of evaluating the right hand side of (\ref{QRflow})\footnote{we
ignore the expense of computing $A(t)$, since this is problem
dependent, and we only focus on the linear algebra expense; further,
recall that we only consider the case of $A$ ``dense''}.
\begin{itemize}
\item[(i)] Runge-Kutta schemes at Gauss points used directly on
(\ref{QRflow}); \cite{DRV}.
\item[(ii)] Projected methods.  There are at least two of these:
(a) {\it bad} and (b) {\it good}.  
\begin{itemize}
\item[(a)]  These methods use the original 
differential equation for $X$, integrate this, and then form its
$QR$ factorization at each step.  They are numerically unsound,
since integration of $X$ is quite often an unstable procedure (as when
the problem is exponentially dichotomic).
\item[(b)]  These consist in using any scheme
to integrate (\ref{QRflow}), and then projecting the solution onto
the manifold of the orthonormal matrices.  For example, the projection step 
can be done with a modified Gram-Schmidt procedure (as in \cite{DRV})
or by using the orthogonal polar factor (as in \cite{Higham:TimeStepPresOrth}).  
\end{itemize}
\item[(iii)] Transformation methods.  Here, one transforms the equation
for $Q$, which is an element of the (Lie) group of the orthogonal matrices,
into an equation for a skew-symmetric matrix $Y$, an element of the
underlying (Lie) algebra.  The equation for $Y$ now evolves in a linear
space, and standard discretizations will deliver a skew-symmetric 
approximation: transforming this back to the group gives the desired
approximation for $Q$.  There have been at least two ways in which
this design has been carried out: 
\begin{itemize}
\item[(a)] using near the identity transformations 
and the matrix exponential to get back on the group
(e.g., see \cite{MuKa:RKLie}), and 
\item[(b)] using the Cayley transform
to get to the algebra and then back to the group
(e.g., see \cite{Diele:NumMetUniDiffSys}).
\end{itemize}
\end{itemize}
\item[[$p<n$]]  Here, the solution belongs to the manifold of 
rectangular orthonormal matrices, the so-called {\sl Stiefel manifold}.
We must appreciate that often one has $p\ll n$, and thus we must 
insist that a well conceived method should
have cost proportional to $O(n^2p)$ per step, that is the cost of 
evaluating the right hand side of (\ref{QRflow}).  As it 
turns out, and to the best of our understanding,
this restriction rules out most of the methods we listed above for the
$p=n$ case, with the exception of projection methods.  However, the 
reasons why the other methods must be discarded are different from one 
another: Gauss RK methods applied directly to (\ref{QRflow})
do not maintain orthogonality 
in the case $p<n$ (see \cite{DiVa:CompLyaExp}), and
transformation methods require going from the group structure to
the algebra and back: apparently, this requires an
$O(n^3)$ expense at some level.  Notice that Gauss RK schemes on 
(\ref{QR-BR}) do maintain orthogonality, but it is not clear to
us that they can be implemented so to require a $O(n^2p)$ expense
per step without forcing a possibly severe stepsize restriction;
for this reason in \cite{BridReich:Stiefel} the authors adopt a 
stabilization procedure which allows for explicit schemes to be
used on (\ref{QR-BR}); in essence, their method becomes akin to
inexact projection methods.
\end{itemize}
\medskip

Based upon the above discussion, and with our present knowledge and
understanding, it would seem that projection methods
are the only survivors amongst the methods recalled above.  
We will see below that the techniques put
forward in \cite{DiVa:ONFact} and revisited 
in this paper are a better alternative,
but before doing so, let us point out some inherent characteristics
and potential limitations of projection methods.  
\begin{itemize}
\item[(i)]  Perhaps the major obstacle to use the ({\sl good}) 
projection methods is of theoretical nature:
projection methods are hard to analyze.
Of course, if any method of order $s$, say, is used 
for (\ref{QRflow}), then the projected solution is also of order $s$.  
However, the projection step itself is a discontinuous operator, and
this has been cause for some worries.  
\item[(ii)]  Part of the worries are probably caused by
the difficulty of getting a backward error statement for projection
methods.  In the numerical analysis of differential equations, 
a backward error analysis consists in the realization
that one has solved (at any given order) a problem which is $O(h^s)$
close to the original one.  We refer to 
\cite{HaLu:lifespan} for this
point of view, which has proven valuable for many problems, and most
notably for Hamiltonian systems.  However, this point of view
is perhaps less relevant in the context under examination
here.  After all, suppose we obtain some orthonormal $\tilde Q$ at the 
grid points $t_k$ instead of the exact $Q$, and some $\tilde R$, instead 
of the exact $R$. Naturally, we can assume that both $\tilde Q$ and
$\tilde R$ are $O(h^s)$ approximation to the exact matrices. Therefore, 
we have triangularized some $\tilde Y(t_k)$ instead of $Y(t_k)$, and 
obviously $\tilde Y(t_k)$ is $O(h^s)$ close to $Y(t_k)$.  In general,
should we expect something better?
\item[(iii)]  In adaptive integration mode, there are some
undesirable features of projection methods.  For the sake of clarity,
suppose that the time stepping strategy is based on two formulas
of different orders.   If we control the error on the unprojected
solution, then it may occasionally happen that we will reject a step
which would not have been rejected if we had checked the error 
on the projected values (or, similarly, we may accept a step which
would have not been accepted for the projected values).  On the other
hand, if we control the error on the projected values, then we essentially
increase the work, because now two projections have to be performed.  
Finally, and regardless of how we control the error, it is undesirable 
that when we end up rejecting a step we had to approximate all of $Q$
in the first place.  It would have been preferable to realize that the
step was not going to be successful ahead of computing all of $Q$.  
Admittedly, this may sound as a strange request, but we'll see
below that it is possible to achieve it with well designed methods.
\item[(iv)]  To conclude, and in spite of the above observations,
we must say that our experience with projection methods
has been quite positive.  For this reason, in Section \ref{s4}, we
will compare performance of our new codes with a projection method.
\end{itemize}
\medskip

We are now ready to list which properties --in our opinion-- a method
to solve (\ref{QRflow}) should have.  At this point, the word ``method''
must be read as: ``formulation of the task in a form which in 
principle allows 
for the properties below to be satisfied''.  In so doing, we maintain
the freedom of delaying discussion of which formulas will
be used in practice.

\begin{enumerate}
\item The method must deliver an orthonormal approximation, at the very
least at the mesh points.  If a RK type integrator is
adopted, the method must deliver (at no added cost) an orthonormal 
approximation also at the internal abscissas.  This will guarantee
that forming the matrix $\tilde A$ at the internal abscissas
is a numerically stable procedure.
\item The method must be flexible enough to handle without modifications
both cases $p=n$ and $p<n$.
\item The method must have a cost per step of $O(n^2p)$, and never
require a $O(n^3)$ expense when $p<n$.  
\item The method must be based on numerically sound coordinates.
\item The method must allow for integration of the relevant differential
equations with theoretical order restrictions given only by the
degree of differentiability of $Q$.
\item The method must be well suited also in adaptive step-size mode:
(i) we want to be able to increase the stepsize if $Q$ evolves slowly,
(ii) we want to be able to reject a step which is going to be unsuccessful
as quickly as possible, possibly (much) earlier than having completed 
approximation for all of $Q$.
\item The method must be flexible enough to allow for increasing 
the number of columns of $A$ which we want to triangularize, without
having to restart the entire triangularization process
from the beginning, and it should
be powerful enough to exploit the work already done.
\end{enumerate}
\medskip

As we will see, the methods we laid down in 
\cite{DiVa:ONFact} achieve the above points. These methods are based upon
writing $Q$ as product of Householder or 
Givens transformations, and since, in general, 
these elementary transformations do not vary smoothly on
the whole interval of interest, we have to adopt a 
``reimbedding'' strategy in order to obtain a well defined process. 
We refer to \cite{DiVa:ONFact} for the original derivation of this approach, 
which we will review in Sections \ref{HousTr} and \ref{GivTr}.

Now we consider a different interpretation of these methods
based on Householder or Givens transformations.
In what follows, we will assume that we have to integrate
(\ref{QRflow}) when $p<n$; trivial simplifications take place
if $p=n$.

We begin observing that, in principle, the sought solution $Q$ can 
always be written in the redundant way
$Q(t)=\bmat{Q(t) & Q^\perp(t)}\bmat{I_p \cr 0}=:U(t)\bmat{I_p \cr 0},
\ \forall t$, where $Q^\perp$
is the orthogonal complement of $Q$.  Thus, we can think
--at least in theory-- to have the following representation for
the sought solution $Q$:
\begin{equation}\label{transQ}
Q(t) \ = \ U(t_k,t_0)Q(t,t_k) \ ,\ t\ge t_k\ ,\,\,\, 
Q(t_k,t_k)\ = \  \bmat{I_p \cr 0}\,.
\end{equation}
In (\ref{transQ}), $U(t_k,t_0)\in \R^{n\times n}$ and
$Q(t,t_k)\in \R^{n\times p}$.  The matrix $U(t_k,t_0)$
is at once comprising $Q$, at $t_k$, and its orthogonal complement.
Now, let $U_k:=U(t_k,t_0)$, and set
$$\hat A(t)\ = \ U_k^TA(t)U_k\,,\,\, t\ge t_k\,.$$
With this, we have that the equation satisfied
by $Q(t,t_k)$ is again (\ref{QRflow}) with $\hat A$ 
replacing $A$ there.  With abuse of notation, if we still
call $Q(t)$ what is really $Q(t,t_k)$, for $t\ge t_k$,
we would then have
\begin{equation}\label{QRflow2}
\dot Q = \hat A Q-QQ^T\hat A Q+Q\hat S\,,\,\ 
\hat S_{ij} =\cases{(Q^T\hat AQ)_{ij},& $i> j$ \cr 
                 0,  &  $i=j$ \cr
   -(Q^T\hat AQ)_{ji},& $i<j$}\,\,, 
\ Q(t_k)=\bmat{I_p \cr 0}\,.
\end{equation}
(The equation for $R$, (\ref{Rflow}), is also modified trivially
with $\hat A$ replacing $A$.)

Thus, quite clearly, knowledge of $U_k$ would allow us to integrate 
the equation (\ref{QRflow2}), 
starting ``near the identity''; ostensibly, this can be done
with a number of choices for representing $Q(t,t_k)$ valid
locally, if not globally (recall our discussion on local
charts).  Now, suppose
we can solve (\ref{QRflow2}) from $t_k$ to $t_{k+1}$, obtaining
$Q(t_{k+1})$,
and that at the same time we are able to obtain also
the orthogonal complement of $Q(t_{k+1},t_k)$, $Q^\perp(t_{k+1},t_k)$;
then, we could form the matrix
$U_{k+1}=\bmat{Q(t_{k+1}) & Q^\perp(t_{k+1})}\ U_k$, redefine
$\hat A$ and $Q(t,t_{k+1})$ accordingly for $t\ge t_{k+1}$ 
and then again formulate a differential equation for $Q(t_,t_{k+1})$
identical to (\ref{QRflow2}) for $t\ge t_{k+1}$.  We could then 
repeat this basic setup until we arrive at $t_f$.
Proceeding this way, we would always have to solve differential
equations with initial conditions $\bmat{I_p \cr 0}$.
Naturally, in practice we will only have computed approximations,
and not exact values, but the setup remains unchanged.  

Notice that, unless we explicitly compute
also $Q^\perp$ so to triangularize all of $A$ (and this
would cost $O(n^3)$), the transformed
matrix $\hat A$, at the $t_k$'s, is upper triangular only in its 
left (n,p) part, in particular its bottom (n-p,n-p) block is
not triangular.  But, suppose that we now decide that
we really needed to triangularize more columns of the fundamental
matrix solution $X$ (or, which is the same, wanted to transform
larger part of $A$ to upper triangular).  If, somehow,
we managed to keep track
not only of $Q(t)$, but also of its orthogonal complement, then
we could work only on the bottom (n-p,n-p) piece of the matrix
$\hat A$.  In other words, at the price of added memory of course,
we would avoid having to restart from scratch.  
But, is it possible to obtain at once information on 
$Q(t_{k+1},t_k)$ and $Q^\perp(t_{k+1},t_k)$, for a cost proportional
to $O(n^2p)$ per step?  This is where we can in principle
exploit the representation for $Q(t,t_k)$ in terms of
either Householder or Givens transformations.  In fact, with these
choices, one does get information on both $Q(t,t_k)$ and its 
orthogonal complement at the same time, at the expense 
of computing only $Q(t_{k+1},t_k)$.  

\smallskip
\noindent{\bf Warning}.
It should be stressed right away that Householder (or Givens) transformations 
are not globally defined; this means that if we want to represent $Q(t)$ 
in (\ref{transQ}) by using these transformations, in general we cannot expect 
to have a globally smooth representation.
As we showed in \cite{DiVa:ONFact}, it is a trivial matter to recover
for free a smooth representation for $Q(t)$.  However, it is {\bf not 
trivial} at all to recover for free a smooth representation for 
$Q^\perp(t)$ (see also \cite{ColemanSorensen:ONBasisForNullSpace}).  After all,
it is hardly imaginable that one can obtain a smooth representation for 
the orthogonal complement at the price of only getting $Q$!  Regardless,
lack of smoothness in the obtainable representation of $Q^\perp$
is not a concern in the contexts of which we are aware.  For example,
if the orthogonal complement is needed because we intend to carry further 
the triangularization process of $A$, then this can be done without restarting
the entire computation from the beginning (but it is not a trivial
implementation to do, since we need to save in memory all quantities which
have been computed thus far).
\bigskip

Now, since Householder and Givens transformations are not unique
(there is a sign ambiguity for each Householder reflection, and
a decision to be made on the order in which we apply the
Givens rotations), we must specify the way in which this
lack of uniqueness is resolved.  This extra freedom means 
that a choice must be made between different coordinates systems
to represent the same object.  The
best way to resolve this ambiguity is to make sure that we
pick up the soundest coordinates from the numerical point of view.

These ideas have been carried out (in a different
notation) in \cite{DiVa:ONFact}.  There, we wrote the solution
$Q(t)$ of (\ref{QRflow}), locally, as product of elementary
Householder or Givens transformations.  We are ready to
review the 
basic algorithms based on Householder and Givens transformations.
We will not explicitly take
advantage of the representation (\ref{transQ}), but, by
writing $Q(t)$ solution of (\ref{QRflow})
as product of elementary orthogonal matrices of Householder
or Givens type, we are conceptually representing
$Q$ as in (\ref{transQ}), with nonsmooth $U$.

\subsection{Householder transformations}\label{HousTr}

Suppose we are at $t_k$ and that we know $X(t_k)$
(e.g., $t_k=t_0$).  Then,
to find $Q(t)$ such that $Q^T(t)X(t)=R(t)$, for $t\ge t_k$, we look for
$Q^T(t)=P_p(t)\cdots P_1(t)$, with $P_i(t)=P_i^T(t)$ 
the Householder matrices 
$$P_i(t)=\bmat{I_{i-1} & 0 \cr 0 & Q_i(t) \cr }\,,\,\,\,
Q_i(t)=I-\frac{2}{u_i^T(t) u_i(t)}u_i(t)u_i(t)^T \,.$$
After $(i-1)$ transformations, the matrix $X$ got transformed into 
$P_{i-1}\dots P_1X$
and its first $(i-1)$ columns have been triangularized.  Let us still
call $X$ the transformed matrix, and let $x_i=X(i:n,i)$ be its 
$i$-th column we need to triangularize.  This is the role of $P_i$.
So, we will set 
\begin{equation}\label{u-variable}
u_i(t)=x_i(t)-\sigma_i \|x_i\|e_1\,,  
\end{equation}
and continue the triangularization process.  The standard textbook
choice for the value of $\sigma_i$
is (see \cite{GolubVanLoan:MatComp}):
\begin{equation}\label{idealsigma}
\sigma_i:=\cases{-1, &if $e_1^Tx_i(t_k)\ge 0$, \cr
1, &if $e_1^Tx_i(t_k)<0$. \cr} 
\end{equation}
This is the idea.  But, of course, we do not have $X$ except at $t_0$.  
In \cite{DiVa:ONFact}
we showed that differential equations can be set for the $u_i$'s directly,
and for the norm of the vectors $x_i$,
requiring only knowledge of the coefficient matrix $A$.  Moreover,
we also derived differential equations for the Householder transformations
in different variables, namely for
\begin{equation}\label{v-variable}
v_i := \frac{u_i}{\norm{u_i}} \,,
\qquad Q_i=I-2v_iv_i^T\,,
\end{equation}
which are better scaled variables, since $v_i^Tv_i=1$.   

We now recall the differential equations for the $u$ and $v$-variables.
For simplicity, we omit the subindices, and thus use the
notation $u$ for $u_i$, etc., and also use $A$ for
$A(i:n,i:n)$, where the matrix $A$ has been progressively
modified by the accumulated transformations:
\begin{equation}\label{AP-upd}
(A,P_j){\rm -update}\,:\,\,\, 
A(t):=P_j(t)A(t)P_j(t)-P_j(t)\dot P_j(t)\,,\,\, j=1,\dots, i-1\,. 
\end{equation}
With this understanding, for the $u$-variables we have:
\begin{equation}\label{u-eq}
\frac{d}{dt}\bmat{u \cr \sigma \|x\| \cr} =
\bmat{A-2e_1e_1^T A_s & (A-e_1 e_1^T A_s)e_1 \cr 2e_1^T 
A_s & e_1^T A_se_1 \cr } \, \bmat {u \cr \sigma \|x\| \cr} 
-{{u^TA_su}\over {\sigma \|x\|}}\bmat{e_1 \cr -1 \cr } \,. 
\end{equation}
Here, $A_s=1/2(A+A^T)$, $e_1$ is the first unit vector (of
appropriate dimension), and (\ref{u-eq}) must be understood as
a differential equation for the $u_i,\ i=1,\dots,p$.

For the $v$-variables, we have the following.
Partition $v$ as $v=\bmat{(e_1^Tv) \cr \hat v}$, let 
$\bmat{a_{11}\cr \hat a_1}$ be the first column of $A$,
$(0,\hat a_{1,a}^T)$ be the first row of ${1\over 2}(A-A^T)$, 
$\bmat{a_{11}\cr \hat a_{1,s}}$ be the first column of $A_s$,
and let $\hat A$ and $\hat A_s$ be the submatrices obtained from
$A$ and $A_s$, respectively, by deleting the first row and column.
Then, we have
\begin{equation}\label{v-eq}
\frac{d}{dt}\bmat{(e_1^Tv) \cr \hat v}=\bmat { 0 & c^T-b^T \cr
b-c & S-S^T \cr} \, \bmat{(e_1^Tv)\cr \hat v}\,, 
\end{equation}
where we have set
\begin{eqnarray*}
&b=\frac{2(e_1^Tv)^2-1}{2}\hat a_1+\hat A\hat v (e_1^Tv)\,,\,\,
c^T:=(-\hat a_{1,a}^T\hat v+\alpha )\hat v^T\,,\,\,
S=(\frac{2(e_1^Tv)^2-1}{2(e_1^Tv)} \hat a_1 + \hat A \hat v)\hat v^T\,,\cr
&{\rm and}\,\,\,
\alpha:=(a_{11}(e_1^Tv)+\hat a_{1,s}^T\hat v)(2(e_1^Tv)^2-1)+2(e_1^Tv)
\hat v^T(\hat a_{1,s}(e_1^Tv)+ \hat A_s\hat v)\,.\end{eqnarray*}

The differential equations (\ref{u-eq}) and (\ref{v-eq}) can easily
be supplied with initial conditions at $t_0$ since $X_0$ is known, and we 
can find $Q_0$ via Householder transformations (in either $u$ or $v$
formulation, with the $\sigma$'s satisfying (\ref{idealsigma})). However,
to describe the typical step between $t_k$ and $t_{k+1}=t_k+h_k$, and
regardless of the choice adopted between $u$ or $v$-variables, the
expression (\ref{idealsigma}) used for choosing the
sign of $\sigma$ must be modified since in
practice we do not have $X(t_k)$.  The way it is done is to enforce
(\ref{idealsigma}), but by only keeping track of the 
transformations.  It goes as follows.

Suppose we have found the Householder matrices at $t_k$, coming from
$t_{k-1}$, call them $P_i^{(k-1)}(t_k)$.  Call
$P^{(k)}_i(t_k)=\bmat{I_{i-1} & 0 \cr 0 & Q^{(k)}_i(t_k)\cr}$
the possibly different initial condition for the Householder
matrices (that is, different $u$'s or $v$'s and $\sigma$'s)
we need in order to step past $t_k$.  Let $K^{(k)}_0=I_n$. 
Inductively define $\sigma^{(k)}_i\,,\,\,i=1,\dots,p$, as follows:
\begin{equation}\label{sigmas}
\sigma^{(k)}_i:=\cases{-1, &if 
$\sigma^{(k-1)}_ie_1^TK^{(k)}_{i-1}Q^{(k-1)}_ie_1\ge 0\,$, \cr
+1, &if $\sigma^{(k-1)}_ie_1^TK^{(k)}_{i-1}Q^{(k-1)}_ie_1<0\,$, \cr}\, 
\end{equation}
where the matrices 
$K^{(k)}_{i-1}\in \R^{n-i+1,n-i+1}\,,\,\,i=2,\dots,p$, are defined by
$$Q^{(k)}_{i-1}(t_k)K^{(k)}_{i-2}Q^{(k-1)}_{i-1}(t_k)=\bmat{
\sigma^{(k-1)}_{i-1} \sigma^{(k)}_{i-1} & 0 \cr
0 & K^{(k)}_{i-1} \cr} \, ,   $$
and, for $i=1,\dots, p$, 
the matrix $Q^{(k)}_i(t_k)$ is obtained so that
$$K^{(k)}_{i-1}Q_i^{(k-1)}(t_k)\sigma^{(k-1)}_ie_1=
Q_i^{(k)}(t_k)\sigma^{(k)}_ie_1\,. $$
Thus, we need to find a\footnote{not unique, there is typo
in \cite{DiVa:ONFact}} Householder transformation which
transforms the left-hand side of this last equation into $\sigma^{(k)}_ie_1$.  
This trivially gives new ICs for the $v_i^{(k)}(t_k)$; the sign ambiguity
in the vector $v_i^{(k)}(t_k)$ is resolved by forcing the sign to
that of the first component of $v_i^{(k-1)}(t_k)$.  
Thus, we can prescribe new ICs for the $v_i^{(k)}(t_k)$, and from these
it is simple to write ICs for the $u_i$'s:
$$u_i^{(k)}(t_k)=-2\sigma^{(k)}_i\|x_i(t_k)\|(e_1^Tv_i^{(k)}(t_k))
v^{(k)}_i(t_k)\,.$$

The proof of the following Lemma is omitted since it amounts to a 
simple verification.
\begin{lem}\label{bestcoord1}
The choice (\ref{sigmas}) is the same as (\ref{idealsigma}).
\end{lem}
\medskip

\begin{rem}\label{bestcoordH}\rm
Let us substantiate the claim that the choice (\ref{idealsigma}) ensures
that we are using sound coordinates from a numerical point of 
view.  In linear algebra, see \cite{GolubVanLoan:MatComp}, the 
choice of signs as in (\ref{idealsigma}) is justified in order to
avoid subtraction of (possibly) nearly equal numbers.  Of course,
this is still true in our context, since we will need to form the
reflectors.  But there is also another aspect to take into account
in the present context.
We restrict to the $v$-coordinates, since these must be generally 
preferred to the $u$-coordinates (but see the $w$-coordinates
of (\ref{w-variable}) in Section \ref{s3}).  We seek an orthonormal 
function $Q$, and we are choosing local coordinates to represent it.  
Obviously, every column of $Q$, $Qe_i$, has unit length, and we represent 
these columns as $P_i\dots P_1e_i$; to do this, we find the $v_i$'s,
each of which has itself unit length, by integrating (\ref{v-eq}).  From
the differential equations (\ref{v-eq}), we observe that there is
a division by $e_1^Tv_i$ in forming the vector
$\frac{2(e_1^Tv_i)^2-1}{2e_1^Tv_i}\hat v$ of $S$ (here,
$e_1=(1,0,\dots,0)^T\in\R^{n-i+1}$).   For stability, we must ensure to 
avoid division by small numbers.   But, with some simple algebra,
one can see that: \hfill\break
\noindent \sl The choice (\ref{sigmas}) is equivalent to having
\begin{equation}\label{wellscaledH}
(e_1^Tv_i)^2\ge \sum_{j=2}^{n-i+1}(e_j^Tv_i)^2\ ,\,\,
i=1,\dots, p\,.
\end{equation}
\rm In particular, the vector $\frac{2(e_1^Tv_i)^2-1}{2e_1^Tv_i}\hat v$ is 
as well scaled as generally possible.
Moreover, (\ref{wellscaledH}) can be used to decide if the current 
Householder frames are numerically stable or not.

In the $u$-variables, it is easier to check (\ref{idealsigma}) directly,
since $e_1^Tx_i=e_1^Tu_i+\sigma_i \|x_i\|$; thus, as long as 
\begin{equation}\label{wellscaledHu}
\sigma_i(e_1^Tu_i+\sigma_i \|x_i\|) < 0 \ ,\,\,
i=1,\dots, p\,,
\end{equation}
there is no need to change the current frame.
\end{rem}
\bigskip

We now summarize the overall strategy on a typical step from 
$t_k$ to $t_{k+1}=t_k+h_k$.  In the next section, we will pay 
closer attention to all aspects of this strategy.

\smallskip
\centerline{{\bf Householder on $[t_k, t_{k+1}]$. }\label{HousAlg}}
\smallskip
\noindent {\tt INPUT}: $t_k$, $h_k>0$, initial conditions 
$P^{(k-1)}_i(t_k)\,,\,\, i=1,\dots, p$ (i.e., the vectors
$u^{(k-1)}_i$ and $\|x_i^{(k-1)}\|$ or the vectors $v^{(k-1)}_i$ at $t_k$),
and the $\sigma^{(k-1)}_i$.  
\begin{itemize}
\item[(1)] For $i=1,\dots, p$, 
check to see if (\ref{wellscaledH}) (or (\ref{wellscaledHu})) holds true.
If it fails, redefine $\sigma^{(k)}_i$ according to (\ref{sigmas}) and 
determine new $P_i^{(k)}(t_k)=\bmat{I_{i-1} & 0 \cr
0 & Q^{(k)}_l(t_k) \cr}$ accordingly (redefine
$u_i^{(k)}(t_k)$ or $v_i^{(k)}(t_k)$), for $i=1,\dots, p$
\item{} For $i=1,\dots, p$
\begin{itemize}
\item[(2)] Let $A=A(i:n,i:n)$
\item[(3)] Find the Householder transformation $P_i^{(k)}(t)$
by integrating (\ref{u-eq}) or (\ref{v-eq}) on $[t_k,t_{k+1}]$ 
\item[(4)] Do an $(A,P_i)$ update (\ref{AP-upd})
\end{itemize}
\item{} Endfor $i$.
\end{itemize}

\noindent {\tt OUTPUT}: $Q^{(k)}(t_{k+1})^T=P^{(k)}_p(t_{k+1})\cdots 
P^{(k)}_1(t_{k+1})$, is such that
$Q^{(k)}(t_{k+1})^TX(t_{k+1})$ is triangular.

\bigskip

\subsection{Givens transformations}\label{GivTr}
Suppose at $t_k$ we know $X(t_k)$
(e.g., $t_k=t_0$).  
To find $Q(t)$ such that $Q^T(t)X(t)=R(t)$, for $t\ge t_k$, we look for
$Q(t)=Q_1(t)\cdots Q_p(t)$, where $Q_i(t)$ is of the form
$Q_i(t)=\bmat{I_{i-1} & 0 \cr 0 & G_i(t)}$, and each $G_i$,
$i=1,\dots, p$, is the product of elementary
planar rotations (Givens, or Jacobi, transformations) of the type
$$Q_{ij}(t)=I-(e_1e_1^T+e_je_j^T)+G_{ij}\ ,\,\
G_{ij}=c_{ij}(e_1e_1^T+e_je_j^T)-s_{ij}(e_1e_j^T-e_je_1^T)\ ,$$
for $j=2, \dots, n-i+1$.
Above, we have used $c_{ij}$ and $s_{ij}$ to express
$\cos (\theta_{ij}(t))$ and $s_{ij}=\sin(\theta_{ij}(t))$, respectively,
where the function $\theta_{ij}(t)$ needs to be found.  Now, suppose 
we have triangularized the first $i-1$ columns of $X$, and still
call $X$ the transformed matrix.  The role of $G_i$ is to
triangularize $x_i:=X(i:n,i)$, the $i$-th column of the unreduced part 
of $X$.  

In the standard linear algebra setting (see \cite{GolubVanLoan:MatComp}), 
the rotators are safely applied in their natural
sequence $Q_{i,i+1},\dots, Q_{i,n}$.  But this may lead to instabilities
in our time dependent setting!  In fact, the specification of the order 
in which the rotators $Q_{i,i+1},\dots, Q_{i,n}$ are applied turns
out to be of utmost importance for numerical stability, and therefore
for accuracy and efficiency.  We adopted the following strategy.

\noindent{\sl Order of rotators}.  Our strategy is 
\begin{equation}\label{rot-strat}
\, {\rm First \,\, rotator \,\, must \,\,
annihilate \,\, largest \,\, entry \,\, of}\,\, x_i(2:n-i+1)\ .
\end{equation}  
To be precise, define $l$ to be the largest entry in
absolute value of $x_i(2:n-i+1)$:
\begin{equation}\label{indexl}
l\ :\,\ X_{i+l-1,i}=\max_{2\le j\le n-i+1} |X_{i+j-1,i}|\,,
\end{equation}
and define the index array $\pi_i$ as
\begin{equation}\label{index-array}
\pi_i=[1,l,2,\dots,l-1,l+1,\dots,n-i+1]\,.
\end{equation}
Then, define (the ordering of the rotators) $G_i$ as
\begin{equation}\label{order-rot}
G_i \ = \ Q_{i,\pi_i(2)}\cdots Q_{i,\pi_i(n-i+1)}\,.
\end{equation}
\begin{rem}\rm 
There seems to be no need to further refine (\ref{rot-strat}) by
selecting the second rotator to annihilate the largest entry
of the unreduced part, and so forth.  
\end{rem}

In \cite{DiVa:ONFact}, we derived differential equations for the elementary
rotators, that is for the $\theta_{ij}$ or for the corresponding
$(\cos,\sin)$ pairs.  To recall, omitting the row index $i$
(i.e., using $\theta_j$ for $\theta_{ij}$, etc.), these are
\begin{equation}\label{DEtheta}
\bmat{c_{\pi_i(3)}\cdots c_{\pi_i(n-i+1)} \ \frac{d}{dt}\theta_{\pi_i(2)} \cr
c_{\pi_i(4)}\cdots c_{\pi_i(n-i+1)} \ \frac{d}{dt} \theta_{\pi_i(3)} \cr
\vdots  \cr
c_{\pi_i(n-i+1)} \ \frac{d}{dt} \theta_{\pi_i(n-i)}  \cr
\frac{d}{dt} \theta_{\pi_i(n-i+1)} \cr}=\ \bmat{\alpha_{\pi_i(2)} \cr 
\alpha_{\pi_i(3)} \cr \vdots \cr \alpha_{\pi_i(n-i)} \cr
\alpha_{\pi_i(n-i+1)} \cr}\,.
\end{equation}
Here, for $j=2,\dots,n-i+1$, we have set
$$\alpha_{\pi_i(j)}(t)=e_{\pi_i(j)}^T \bigl[Q_{i,\pi_i(n-i+1)}^T\cdots 
Q_{i,\pi_i(2)}^T A
Q_{i,\pi_i(2)}\cdots Q_{i,\pi_i(n-i+1)}\bigr]e_1\,,$$
and $A$ is really $A(i:n,i:n)$, which has been progressively modified
by the accumulated transformations:
\begin{equation}\label{AQ-upd}
(A,Q_j){\rm -update}\,:\,\,
A(t):=Q_j(t)^TA(t)Q_j(t)-Q_j^T(t)\dot Q_j(t)\ , \ \, j=1,\dots, i-1\ . 
\end{equation}
Of course, since $\frac{d}{dt}\bmat{\cos(\theta_{ij}) \cr
\sin(\theta_{ij})}=\bmat{-\sin(\theta_{ij}) \cr \cos(\theta_{ij})}
\frac{d}{dt}\theta_{ij}$, from (\ref{DEtheta}) it is trivial to write
differential equations for the $(\cos,\sin)$ pairs directly:
\begin{equation}\label{DEcs}
\bmat{c_{\pi_i(3)}\cdots c_{\pi_i(n-i+1)} \frac{d}{dt}\bmat{c_{\pi_i(2)} \cr 
s_{\pi_i(2)}}\cr
  \vdots \cr   \cr
c_{\pi_i(n-i+1)} \frac{d}{dt}\bmat{c_{\pi_i(n-i)} \cr s_{\pi_i(n-i)} \cr}\cr
\frac{d}{dt}\bmat{c_{\pi_i(n-i+1)} \cr s_{\pi_i(n-i+1)}}}=\bmat{
\bmat{0& - \alpha_{\pi_i(2)} \cr \alpha_{\pi_i(2)} & 0} 
\bmat{c_{\pi_i(2)} \cr s_{\pi_i(2)}} 
\cr  \cr 
\vdots \cr   \cr 
\bmat{0 & -\alpha_{\pi_i(n-i)} \cr \alpha_{\pi_i(n-i)} & 0 \cr} 
\bmat{c_{\pi_i(n-i)} \cr  s_{\pi_i(n-i)}}  \cr 
\bmat{0 & -\alpha_{\pi_i(n-i+1)} \cr
\alpha_{\pi_i(n-i+1)} & 0 \cr} \bmat{c_{\pi_i(n-i+1)} \cr  s_{\pi_i(n-i+1)}}} . 
\end{equation}

To supply initial conditions at $t_0$ for the
differential equations (\ref{DEtheta}) and/or (\ref{DEcs}), we enforce
(\ref{rot-strat}); that is, use
(\ref{indexl}), (\ref{index-array}), and apply the rotators in the order
specified by (\ref{order-rot}).  At each application of an elementary
rotator on $x_i$ at $t_0$, there is a sign ambiguity;
we resolve this ambiguity by forcing the first entry of
$x_i$ to be always positive as the rotators are applied
(and thus, equal to $\|x_i\|$ after all the
elementary rotators $Q_{i,\pi_i(2)},\dots, Q_{i,\pi_i(n-i+1)}$,
are applied).  This way, we can start the integration.

But, to describe the typical step between $t_k$ and $t_{k+1}=t_k+h_k$, 
the strategy just outlined, in particular the
expressions (\ref{indexl}) and (\ref{index-array}), must be modified, 
since we do not have $X$ at $t_k$.  To understand the way we do this,
we first need the following remark.

\begin{rem}\label{bestcoordG}\rm
Givens transformation provide another
example of the situation we described at the beginning of this work, when
we discussed how a trajectory on a smooth manifold
can be parametrized by overlapping
local charts.  Then, we said that an appropriate way to choose local 
coordinates was to enforce numerical stability.  We now justify 
how this can be achieved with the choices we have adopted.
Without loss of generality, assume that $\pi_i=[1,2,\dots, n-i+1]$ (if
not, relabel the indices accordingly).
From the differential equations (\ref{DEtheta}) and (\ref{DEcs}), we 
notice that multiplication by the inverse of the diagonal matrix
$$\diag\bigl(c_{3}\cdots c_{n-i+1}, \ c_{4}\cdots c_{n-i+1}, \ \dots, \
c_{n-i+1},\ 1\bigr)$$ 
is taking place.  
For numerical stability, we need to avoid division by small numbers.
Clearly, the smallest number by which we are dividing (in absolute sense) 
is $c_{3}\cdots c_{n-i+1}$.  
Let $x_{i,j}$ denote the $j$th component of $x_i$, $j=1,2,\dots, n-i+1$.
Then, using rotators to triangularize $x_i$, we have
$$c_j^2=\frac{(\prod_{l=1}^{j-1} x_{i,l}^2)}
           {(\prod_{l=1}^j x_{i,l}^2)}\,,$$
and so
$$c_{3}^2\cdots c_{n-i+1}^2 = \frac{(x_{i,1}^2+x_{i,2}^2)}
           {\norm{x_i}^2}\,.$$
Now, as long as 
\begin{equation}\label{x-ineq}
x_{i,1}^2+x_{i,2}^2 \ \ge x_{i,j}^2,\,\,j=3,\dots, n-i+1\ ,\,\, 
i=1,\dots, p\, ,
\end{equation}
no division by small number is taking place, since when (\ref{x-ineq})
is satisfied we have
$$c_{3}^2\cdots c_{n-i+1}^2 \ge \frac{1}{(n-i)}\,,$$
which is as nicely bounded away from $0$ as any product of cosines
for Givens' transformations can generally be.
Furthermore, with some simple algebra, one can see that: \hfill\break
\noindent \sl The inequality (\ref{x-ineq}) is equivalent to having
\begin{equation}\label{cheapcos}
\prod_{j=3}^k c_j^2\ \ge \ s_k^2\ ,\,\ k=3, \dots, n-i+1\ ,\,\ 
i=1,\dots, p\,.
\end{equation}
\rm 
Clearly, (\ref{cheapcos}) can be used to decide if the current 
ordering of Givens' transformations is numerically stable or not,
without knowledge of $X$.
\end{rem}
\bigskip

We are ready to describe how we modify the strategy which led us to
(\ref{indexl}) and (\ref{index-array}) after the first step.
First of all, as long as (\ref{cheapcos}) holds, we do not 
change the present ordering of the rotators.  In case (\ref{cheapcos}) 
fails, we enforce (\ref{rot-strat}), but by only keeping track of the 
transformations, in the following way.
Suppose we have found the Givens transformation matrices at $t_k$, coming 
from $t_{k-1}$, call them $Q_i^{(k-1)}(t_k)$.  Call
$Q^{(k)}_i(t_k)=\bmat{I_{i-1} & 0 \cr 0 & G^{(k)}_i(t_k)\cr}$
the possibly different initial condition for the new Givens matrices
(that is, different $\theta$'s or cosines/sines)
we need in order to step past $t_k$.
Define $K_0^{(k)}=I_n$, and inductively
define 
\begin{equation}\label{w-i}
w_i:=K_{i-1}^{(k)}G_i^{(k-1)}(t_k) e_1\,,
\,\,\,{\rm for}\,\, i=1,\dots, p \,.
\end{equation}
Let $l$ be the index of the largest entry (in absolute value) of 
$w_i(2:n-i+1)$.  Accordingly, define $\pi_i$ as in (\ref{index-array})
and the ordering for the rotators relative to $G_i^{(k)}(t_k)$ as in
(\ref{order-rot}).  Find initial conditions for the $Q_{i,j}^{(k)}(t_k)$
by enforcing that all vectors 
$\prod_{j=2}^{n-m+1}(Q_{i,\pi_i(j)}^{(k)}(t_k))^Tw_i,$ $m=n-1,\dots, i$, 
have positive first
component.  This way, we define $G_i^{(k)}(t_k)$, hence $Q_i^{(k)}(t_k)$.
Finally, for $i=2,\dots, p$, we define
$K_{i-1}^{(k)}\in \R^{n-i+1, n-i+1}$ from
$$Q_{i-1}^{(k)}(t_k)^T\,
\pmatrix{I_{i-2}& 0 \cr 0 & K_{i-2}^{(k)} \cr } 
Q_{i-1}^{(k-1)}(t_k)=
\pmatrix{I_{i-1}& 0 \cr 0 & K_{i-1}^{(k)} \cr }\,. $$

Again we omit the proof of the following Lemma, since it is a
direct verification.
\begin{lem}\label{bestcoord2}
The choice just described for providing initial conditions on the rotators
at $t_k$ is equivalent to the strategy based on (\ref{order-rot}) and
first paragraph after (\ref{DEcs}).
\end{lem}
\medskip

We summarize the overall strategy on a typical step from 
$t_k$ to $t_{k+1}=t_k+h_k$.  We assume that integration has
been successful up to $t_k$, and that the step $h_k$ has been
selected.

\smallskip
\centerline{{\bf Givens on $[t_k, t_{k+1}]$. }}
\smallskip
\noindent {\tt INPUT}: $t_k$, $h_k>0$, initial conditions 
$Q^{(k-1)}_i(t_k)\,,\,\, i=1,\dots, p$ (i.e., either the 
$(\cos, \sin)$ pairs or the $\theta$ values, and the ordering in 
which they had been applied).
\begin{itemize}
\item[(1)] For $i=1,\dots, p$, 
check to see if (\ref{cheapcos}) holds true.
If it fails, redefine 
the index array $\pi_i$ and initial conditions for
the rotators at $t_k$.  That is, for $i=1,\dots, p$,
let $Q_i^{(k)}=Q_{i,\pi_i(2)}\cdots Q_{i,\pi_i(n+i-1)}$, and find initial
conditions for $Q_i^{(k)}(t_k)$ by bringing $w_i(i:n)$ into $e_1$ so that 
all rotated vectors  
$\prod_{j=2}^{n-m+1}(Q_{i,\pi_i(j)}^{(k)}(t_k))^Tw_i$, 
$m=n-1,\dots, i$, have positive first component.
\item{} For $i=1,\dots, p$
\begin{itemize}
\item[(2)] Let $A=A(i:n,i:n)$.
\item[(3)] Find the transformation $Q_i^{(k)}(t)$
by integrating the differential equations (\ref{DEtheta}) or (\ref{DEcs})
on $[t_k,t_{k+1}]$.
\item[(4)] Do an $(A,Q_i)$ update (\ref{AQ-upd}).
\end{itemize}
\item{} Endfor $i$.
\end{itemize}

\noindent {\tt OUTPUT}: $Q^{(k)}(t_{k+1})=Q^{(k)}_1(t_{k+1})\cdots 
Q^{(k)}_p(t_{k+1})$, is such that
$(Q^{(k)}(t_{k+1}))^TX(t_{k+1})$ is triangular with positive diagonal
entries.

\bigskip
\begin{rem}\rm 
We observe that --even through a change of initial conditions--
with our strategy the signs on the diagonal of $R$, see (\ref{Rflow}),
remain fixed when using rotators to find $Q$.
\end{rem}

\section{Implementation}\label{s3}

Here we describe how we implemented the algorithm put forward
in the previous section.  Before doing so, however, we want to
derive a different formulation for the
Householder transformations, which present some distinct advantages over
both the $u$ and $v$ formulations.  The motivation comes from Remark 
\ref{bestcoordH}.  

We experimented with two different rewritings of the Householder
vectors.
\begin{itemize}
\item[(a)]
With the notation of (\ref{v-eq}), we observe that
$(e_1^Tv)=\sigma(1-\hat v^T\hat v)^{1/2}$.  Thus, one can 
write differential equations just for $\hat v$, and then recover
all of $v$ by using $(e_1^Tv)=\sigma(1-\hat v^T\hat v)^{1/2}$.  On
the surface, this seems to present some advantages,
since there is one less differential equation to solve.
However, there is a potential loss of precision which occurs 
with this approach when forming back $e_1^Tv$.   After all, we 
know that the first component of $v$ dominates all others; thus,
we expect that some of the other components will be small
(and we often observed this to be true in practice).  Indeed, in our
experiments, this formulation gave consistently less precise
results than the $v$ formulation or the $w$ formulation to
be introduced next.  For this reason, henceforth we will not present 
computational results obtained with this formulation.
\item[(b)]  
Consider the new variable $w$:
\begin{equation}\label{w-variable}
w\ = \ \frac{1}{e_1^Tv}\ v\,,\,\, {\rm thus}\,\, \, w\ = \
\bmat{1 \cr \hat w}\,.
\end{equation}
Observe that
$$(e_1^Tv)\ = \ \frac{-\sigma}{\norm{w}}\,,$$
so that the discussion relative to the choice of $\sigma$'s 
(see (\ref{sigmas})) stays unchanged.  Since
$\dot w=\frac{dw}{dt}\ = \ \bmat{0 \cr \frac{d\hat w}{dt}}$, we have one
less differential equation to solve and a simplified form
for the $(A,P_i)$-updates.  Omitting the indices for simplicity,
the form of the typical update becomes 
\begin{eqnarray}\label{w-update}
& PAP-P\dot P\ = \ \\
&= \ A-\frac{2}{w^Tw}( w(w^TA)+(Aw)w^T)+
4\frac{w^TAw}{(w^Tw)^2} ww^T-\frac{2}{w^Tw}
(w\dot w^T-\dot ww^T)\ . \nonumber
\end{eqnarray}
From (\ref{w-update}), it is easy to derive the differential equation
satisfied by $\hat w$.  In fact, this can be obtained from the 
requirement that $(PAP-P\dot P)e_1$ has only its first entry not $0$.
Using the form of $\dot w$, and the same notation used to derive
(\ref{v-eq}), this requirement becomes
$$\hat a_1-\frac{2}{w^Tw}((a_{11}+\hat w^T\hat a_1)\hat w +\hat a_1
+\hat A \hat w)+4\frac{w^TAw}{(w^Tw)^2}\hat w+\frac{2}{w^Tw}
\dot {\hat w} = 0 \,, $$
from which we get the equation for $\hat w$:
\begin{equation}\label{w-eq}
\frac{d \hat w}{dt}\ = \ 
\bigl[a_{11}+\hat w^T\hat a_1 -2\frac{w^TAw}{w^Tw}\bigr] \hat w
+(1-\frac{w^Tw}{2})\hat a_1 +\hat A \hat w \,.
\end{equation}
\end{itemize}
\begin{rem}\label{wscaled}\rm
Notice that the division by $w^Tw$ in (\ref{w-eq}) is perfectly
safe, given the form of $w$.  
Of course, if one uses the $w$-variables, obvious modifications
take place in the skeleton of the algorithm on page \pageref{HousAlg}.
E.g., (\ref{wellscaledH}) now would read
\begin{equation}\label{wellscaledHw}
1-(w_{i+1,i}^2+\dots+w_{n,i}^2)\ \ge \ 0 \,. 
\end{equation}
\end{rem}
\medskip

To sum up, for Householder methods, we have considered three 
possibilities:
(i) $u$-variables, (ii) $v$-variables, and (iii) $w$-variables.  
In the $v$-variables, we need to integrate (\ref{v-eq}) whose solution
has norm $1$ for all $t$, and we must maintain this property under
discretization, at gridpoints.  
In principle, we could use Gauss RK schemes for this
task; we did this in \cite{DiVa:ONFact}, by using Newton's method to solve
the resulting nonlinear system, but this gives a cost of 
$O(pn^3)$ per step which is more than what we are willing to pay.  
Use of the linearly convergent scheme of \cite{DRV}
forces a severe stepsize restriction, which is clearly
incompatible with adaptive time stepping.   For these reasons, 
for the $v$-variables, in this work we propose integrating (\ref{v-eq}), 
with an explicit scheme followed by a renormalization at each step.  The
resulting method is a ``local projection'' method for the $v$-variables
and has cost of $O(n^2p)$ per step.
As far as the integration for the $u$ or $w$-variables, integration
can be carried out with explicit schemes, with no need of renormalization.
Notice, though, that when forming the Householder transformations, and
hence $Q$, one is in fact normalizing things so to keep them orthogonal;
this is obvious: in the matrix $I-2ww^T/w^Tw$, the
term $ww^T/w^Tw$ is indeed an orthogonal projection.

As far as methods based on Givens transformations, so far, we spent more
time trying to obtain good
codes only for the formulation in terms of the $\theta$-variables,
see (\ref{DEtheta}),
and not for the $(\cos,\sin)$-variables.  There are several reasons
for our choice.  First, a simplicity reason: it is simpler to work
with the $\theta$-variables, and then form the rotators.  Second, if
we work with $(\cos,\sin)$, then we must integrate 
(\ref{DEcs}) maintaining the solution of norm $1$: again, we could use Gauss
RK schemes with Newton's method, or a ``local projection'' technique.  
But, in all cases, (\ref{DEcs}) is twice the dimension of the 
system to be solved with the $\theta$-variables.  Repeated use of
the trigonometric identity $\cos^2 \alpha+\sin^2 \alpha=1$ would allow us
to half the dimensions, of course, but at the price of
added nonlinearities.  
Third, there is an accuracy reason to prefer the $\theta$-variables;
in all our tests, they gave more accurate
results than their $(\cos,\sin)$ counterpart.  An explanation
for this fact is the content of the next Lemma.

\begin{lem}\label{thetacc}
Let $\theta$ be the angle associated to the elementary rotation
$Q(\theta)=\bmat{\cos(\theta) & -\sin(\theta) \cr
\sin(\theta) & \cos(\theta)}$.  Let
$\phi$ be an approximation to $\theta$, and 
let $Q(\phi)$ be the elementary rotator associated to $\phi$.  
Consider the error matrix $E:=Q^T(\theta)(Q(\phi)-Q(\theta))$.  If
$\theta-\phi= \eta$, sufficiently small, 
then the error on the diagonal of $E$ is $O(\eta^2)$:
$$E=\bmat{-\frac{\eta^2}{2}+O(\eta^4) & \eta+O(\eta^3) \cr
-\eta + O(\eta^3) & -\frac{\eta^2}{2}+O(\eta^4)}\,.$$
\end{lem}
\begin{proof}
The proof follows from $Q^T(\theta)Q(\phi)=\bmat{\cos(\theta-\phi) & 
\sin(\theta-\phi) \cr -\sin(\theta-\phi) & \cos(\theta-\phi)}$.
\end{proof}

In particular, this Lemma implies that $\|E\|=|\eta|+O(\eta^2)$.  Instead,
an error equal to $\eta$ on $\cos(\theta)-\cos(\phi)$ and on
$\sin(\theta)-\sin(\phi)$, would have rendered $\|E\|=|\eta| \sqrt{2}+
O(\eta^2)$.  In the general case of $Q\in \R^{n\times p}$, with notation 
similar to Lemma \ref{thetacc}, by
taking into account the general form of the matrix $Q$ as product of
rotators, and using Lemma \ref{thetacc} over and over,
it is simple to realize that an error of order $\eta$ on the 
angles still gives an $O(\eta^2)$ error term  on the diagonal of $E$, whereas
an error $O(\eta)$ for the $\cos$ and $\sin$ would give an $O(\eta)$
error term for all entries of $E$.

\begin{rem}\label{needatan}\rm
Of course, 
the $\theta$-variables are more properly seen as variables on a torus.
Indeed, also to avoid pathological cases in which the $\theta$ values 
would grow unbounded, we always renormalize their values to 
$[-\pi, \pi]$, which we do by computing inverse tangent.
This represents a drawback with respect to
working directly with the $(\cos, \sin)$ pairs.
\end{rem}
\begin{rem}\label{minpara}\rm
Householder methods based on the $w$-variables, and Givens' methods
based on the $\theta$-variables parametrize the sought $Q$ using
exactly $p\frac{2n-p-1}{2}$ parameters: the minimal number required.
\end{rem}

\medskip
\noindent{\sl Linear Algebra Involved}.  An advantage of using 
Householder and Givens transformations to represent $Q$ is that we
can take advantage of the many clever ways in which these transformations
can be manipulated, see \cite{GolubVanLoan:MatComp}, and thus obtain
efficient procedures.  For example, we keep the matrix $Q$ in factored
form, and never form it (except when required for output purposes).
Of course, we never perform matrix-matrix multiplications either, but
we exploit inner product arithmetic for Householder matrices and the 
simple structure of the rotators for Givens matrices.

\medskip
\noindent{\sl Discretization Schemes}.  We have chosen to use 
explicit integrators of RK type as the basic schemes.  We are not 
interested in extremely accurate computations; the range of practical
interest for us is between $10^{-2}$ and $10^{-10}$, and our
schemes have been chosen with this accuracy demands in mind.  If one
needs more accurate computation, then different choices would be more
appropriate. We chose 
formulas of order 4 and of order 5 with an associated 
embedded formula, to be used in variable stepsize mode.
The formulas used are the well known Runge-Kutta
3/8-th rule, a scheme of order 4 with an embedded scheme of order 3, 
and the formulas of Dormand-Prince of order 5 with the embedded 
scheme of order 4.  For convenience, we give these pairs below.
The lower order scheme is the one relative to the weights $\hat b$.

\begin{center}
\begin{tabular}{ c | c c c c c}
  $0$ &$0$ &$0$ &$0$ &$0$ &$0$ \\
  $\frac 13$ &$\frac 13$ &$0$ &$0$ &$0$ &$0$ \\
  $\frac 23$ &$-\frac 13$ &$1$ &$0$ &$0$ &$0$ \\
  $1$ &$1$ &$-1$ &$1$ &$0$ &$0$ \\
\hline
  $b$ &$\frac 18$ &$\frac 38$ &$\frac 38$ &$\frac 18$ &$0$ \\
\hline
  $1$ &$\frac 18$ &$\frac 38$ &$\frac 38$ &$\frac 18$ &$0$ \\
\hline
  $\hat b$ &$\frac {1}{12}$ &$\frac 12$ 
&$\frac 14$ &$0$ &$\frac 16$ \\
\end{tabular}
\end{center}
\smallskip
\begin{center}
3/8-th Runge-Kutta 4-3 pair
\end{center}

\bigskip

\begin{center}
\begin{tabular}{ c | c c c c c c c}
  $0$ &$0$ &$0$ &$0$ &$0$ &$0$ &$0$ &$0$\\
  $\frac 15$ &$\frac 15$ &$0$ &$0$ &$0$ &$0$ &$0$ &$0$\\
  $\frac {3}{10}$ &$\frac {3}{40}$ &$\frac{9}{40}$ &$0$ &$0$ &$0$ 
&$0$ &$0$\\
  $\frac 45$ &$\frac {44}{45}$ &$-\frac{56}{15}$ 
&$\frac {32}{9}$ &$0$ &$0$ &$0$ &$0$\\
  $\frac 89$ &$\frac {19372}{6561}$ &$-\frac{25360}{2187}$ 
&$\frac {64448}{6561}$ &$-\frac {212}{729}$ &$0$ &$0$ &$0$\\
  $1$ &$\frac {9017}{3168}$ &$-\frac {355}{33}$ &$\frac {46732}{5247}$ 
&$\frac {49}{176}$ &$-\frac {5103}{18656}$ &$0$ &$0$\\
\hline
  $b$ &$\frac {35}{384}$ &$0$ &$\frac {500}{1113}$ 
&$\frac {125}{192}$ &$-\frac {2187}{6784}$ &$\frac {11}{84}$ &$0$\\
\hline
  $1$  &$\frac {35}{384}$ &$0$ &$\frac {500}{1113}$ 
&$\frac {125}{192}$ &$-\frac {2187}{6784}$ &$\frac {11}{84}$ &$0$\\
\hline
  $\hat b$ &$\frac {5179}{57600}$ &$0$ &$\frac {7571}{16695}$ 
&$\frac {393}{640}$ &$-\frac {92097}{339200}$ &$\frac {187}{2100}$ 
&$\frac {1}{40}$\\
\end{tabular}
\end{center}
\smallskip
\begin{center}
Dormand--Prince 5-4 pair.
\end{center}
\medskip

\begin{rem}\rm
With our schemes based on elementary transformations, one can form
orthogonal approximations also
at the internal RK points, not just at the grid points, at negligible
extra cost.  This is a nice fact, both for dense output purposes, and
for forming $\tilde A$ at the internal points.
\end{rem}

\medskip
\noindent{\it Updates}.
Whenever the updates are required, see 
(\ref{AP-upd}) and (\ref{AQ-upd}),
we do not setup the updated matrices by using exact derivatives
(as we did in \cite{DiVa:ONFact}),
since we can (and do) use the values obtained at the Runge-Kutta
stages; this way, the update does not require
extra function evaluations.

\medskip
\noindent{\it Error Control}.
Finally, about error control for variable time-stepping integration.
Most of the criteria we used are standard.  We do mixed absolute/relative
error control with respect to the input value {\tt TOL}, 
in the way explained in 
\cite[II-4]{HairerNorsettWanner:SODEI}, with some slightly different
heuristics.  For example we use a safety factor of $0.8$ (it is {\tt fac}
in \cite{HairerNorsettWanner:SODEI}), we never allow a new stepsize to 
be bigger than four times the current stepsize, and we choose the initial
stepsize to be ${\mathtt {TOL}}^{1/(q+1)}$, where $q$ is the order of
the estimator (i.e., $q=3$ for the embedded 3/8th rule, and $q=4$ for
the embedded Dormand-Prince rule).  But the most relevant
changes to the standard strategies are due to the nature of our methods.
In fact, triangularizing one column of $X$ at the time, as we do, present 
some important advantages.  First of all, we decide on step-size changes
by monitoring the behavior of the error (in the $\infty$-norm) on all 
columns independently, and then take
the most conservative estimate for the next step.  But, most importantly,
there is no need to complete the entire integration step prior to rejecting
the step!  This is a pleasant outcome of the present implementation, since
often (see Section \ref{s4}) a 
step failure occurs ahead of having completed computation of all of $Q$.

\medskip
\noindent{\it Reimbedding}. 
We have used the construction based on (\ref{sigmas}), or
(\ref{w-i}), and following discussion, only if the tests
(\ref{wellscaledHu}), or (\ref{wellscaledH}), or (\ref{wellscaledHw}),
or (\ref{cheapcos}), for the $u$, $v$, $w$, or $\theta$ variables,
respectively, failed.
We have not succeeded in finding less expensive ways to obtain new
initial conditions in case these tests failed.

\bigskip

\section{Codes \& Examples}\label{s4}
We have written {\tt FORTRAN} codes using the previous
ideas.  In the examples, we will refer to the performance of each
of these codes, which differ by which formulas are used, and whether
or not they are implemented in constant or variable stepsize modes.
We use the following first letter convention:
{\tt u}=$u$-variables,
{\tt v}=$v$-variables, {\tt w}=$w$-variables, {\tt t}=$\theta$-variables.
\begin{itemize}
\item{} Fixed stepsize codes.  
\begin{itemize}
\item{} 3/8th rule: {\tt urk38, vrk38, wrk38, trk38}.  Thus, for example,
{\tt wrk38} is a fixed stepsize implementation using the Runge-Kutta
3/8th rule of the Householder method based on the $w$-variables.
\item{} Dormand-Prince rule: {\tt udp5, vdp5, wdp5, tdp5}.
\end{itemize}
\item{} Variable stepsize codes.  The naming convention is as
above, but now the first letter is a {\tt v} to signify ``variable
stepsize''.
\begin{itemize}
\item{} 3/8th pair: {\tt vurk38, vvrk38, vwrk38, vtrk38}.  For example,
{\tt vvrk38} is the variable stepsize implementation of the 3/8th pair
for the Householder method in $v$-variables.
\item{} Dormand-Prince pair: {\tt vudp5, vvdp5, vwdp5, vtdp5}.
\end{itemize}
\end{itemize}
\medskip

We report on several measures of performance of the above codes.
\begin{itemize}
\item[-] {\tt err}: the error between computed and exact $Q$, if the
exact $Q$ is known.
\item[-] {\tt reimb}: the number of reimbeddings needed for a given method;
we increment the counter every time a change of coordinates is
performed.
\item[-] {\tt rejs/first}: the number of total rejections (in
variable stepsize mode) followed by the rejections occurring while
triangularizing the first column.
\item[-] {\tt cpu}: the total {\tt CPU} time needed to complete a given
run normalized to $1$ for the fastest run on the given problem.
\item[-] {\tt nsteps}: the total number of steps taken (in variable
stepsize mode).
\end{itemize}

For the variable stepsize codes, we also include comparison with a
projected integrator, {\tt prk45}, which integrates (\ref{QRflow})
with the well known (and sophisticated) integrator {\tt RKF45} of
{\tt Netlib}, and then uses modified Gram-Schmidt for the projection.
We recall that {\tt RKF45} is a RK solver of order 5/4 whose performance
is comparable with the Dormand-Prince 5/4 pair we adopted here.

\begin{exm}\label{Ex1}\rm
This is a problem chosen because the underlying fundamental solution is both
exponentially dichotomic and fast rotating.  Bad projection methods have
difficulties on this problem.  Also Gauss RK schemes without
Newton method run into serious stepsize restrictions.  
We have the coefficient matrix 
$$A(t) = \bmat{\beta \cos(2\alpha t) & -\alpha+\beta\sin(2\alpha t)\cr
\alpha+\beta\sin(2\alpha t) & -\beta \cos(2\alpha t)\cr}\,,
$$
and we seek $Q$ associated to the QR factorization of $X:\ X'=AX$, 
$X(0)=I$.  The exact solution is $X(t)=\bmat{\cos(\alpha t) & 
\sin(\alpha t)\cr \sin(\alpha t) & -\cos(\alpha t)\cr}\bmat{
e^{\beta t} & 0 \cr 0 & -e^{-\beta t}}$.  We fix $\beta=100$, 
$\alpha=100$, and consider integration on the interval $[0,b]$ with $b=10$.
The problem should cause difficulties to methods based 
on Householder transformations, because of the fast rotation of $Q$.
This does indeed produce several reimbeddings for Householder
methods, but no appreciable deterioration in accuracy.

Observe that for the {\tt t} methods there is only one equation to
integrate:
$$\dot \theta = \alpha -\beta \sin(2\theta(t))\cos(2\alpha t)+
\beta \sin(2\alpha t)\cos(2\theta(t))\ ,\,\,\, \theta(0)=0\,,$$
and this has the exact solution $\theta(t)=\alpha t$.  So, one may 
expect no error while integrating for $\theta$.  However, the reason
why the {\tt t} methods do not recover the exact solution is because
we automatically renormalize angles to $[-\pi,\pi]$ and this causes
roundoff errors to enter in the picture; if we do not renormalize
the angles, then the exact solution is recovered.

Tables 4.1 and 4.2 summarize the results of our numerical experiments.
For the fixed step methods in Table 4.1 the small error for the {\tt t}
methods is notable and this is mirrored by the superior performance
for the 
variable step {\tt t} methods (see Table 4.2). Note that the {\tt u}
methods fail on this problem.  Also, observe how {\tt prk45} is 
considerably more expensive than the competing 5/4 codes ({\tt vvdp5,
vwdp5}).

\bigskip
\centerline{
\vbox{\tabskip=0pt\offinterlineskip
\def\tblrule{\noalign{\hrule}}
\halign{ 
 \strut#&\vrule#\tabskip=1em plus2em
 &\hfil#\hfil&\vrule#    
 &\hfil#\hfil&\vrule#    
 &\hfil#\hfil&\vrule#    
 &\hfil#\hfil&\vrule#    
 \tabskip=0pt\cr         
\tblrule                      
&&\multispan{7}\hfil Table 4.1. Example \ref{Ex1}: fixed stepsize 
$\Delta t=1.E-3$.
\hfil&\cr\tblrule
\tblrule\tblrule
&&\omit\hidewidth {\tt Meth} \hidewidth&&
  \omit\hidewidth {\tt err} \hidewidth&&
  \omit\hidewidth {\tt reimb} \hidewidth&&
  \omit\hidewidth {\tt cpu} \hidewidth&\cr\tblrule
\tblrule\tblrule
&& {\tt tdp5}  && $3.1E-13$ && $0$   && $17.5$ &\cr\tblrule
&& {\tt trk38} && $3.9E-13$ && $0$   && $14.2$ &\cr\tblrule
&& {\tt udp5}  && $-$       && $-$   && $-$    &\cr\tblrule
&& {\tt urk38} && $-$       && $-$   && $-$    &\cr\tblrule
&& {\tt vdp5}  && $2.5E-9$  && $318$ && $24.0$ &\cr\tblrule
&& {\tt vrk38} && $1.6E-6$  && $318$ && $18.9$ &\cr\tblrule
&& {\tt wdp5}  && $3.9E-8$  && $318$ && $18.9$ &\cr\tblrule
&& {\tt wrk38} && $2.4E-6$  && $318$ && $15.5$ &\cr\tblrule
\cr}}}
\par\bigskip

\bigskip
\centerline{
\vbox{\tabskip=0pt\offinterlineskip
\def\tblrule{\noalign{\hrule}}
\halign{ 
 \strut#&\vrule#\tabskip=1em plus2em
 &\hfil#\hfil&\vrule#    
 &\hfil#\hfil&\vrule#    
 &\hfil#\hfil&\vrule#    
 &\hfil#\hfil&\vrule#    
 &\hfil#\hfil&\vrule#    
 &\hfil#\hfil&\vrule#    
 \tabskip=0pt\cr         
\tblrule                      
&&\multispan{11}\hfil Table 4.2. Example \ref{Ex1}: variable stepsize,
{\tt tol}$=1.E-8$.
\hfil&\cr\tblrule
\tblrule\tblrule
&&\omit\hidewidth {\tt Meth} \hidewidth&&
  \omit\hidewidth {\tt err} \hidewidth&&
  \omit\hidewidth {\tt reimb} \hidewidth&&
  \omit\hidewidth {\tt rejs} \hidewidth&&
  \omit\hidewidth {\tt cpu} \hidewidth&&
  \omit\hidewidth {\tt nsteps} \hidewidth&\cr\tblrule
\tblrule\tblrule
&& {\tt prk45}  && $1.4E-8$ && $-$ && $-$  && $31.1$   && $20803$ &\cr\tblrule
&& {\tt vtdp5}  && $4.6E-8$ && $0$   && $162$  && $1.2$   && $599$ &\cr\tblrule
&& {\tt vtrk38} && $2.5E-8$ && $0$   && $158$  && $1$   && $705$ &\cr\tblrule
&& {\tt vudp5}  && $-$        && $-$   && $-$    && $-$      && $-$ &\cr\tblrule
&& {\tt vurk38} && $-$        && $-$   && $-$    && $-$      && $-$ &\cr\tblrule
&& {\tt vvdp5}  && $3.3E-9$ && $318$ && $0$    && $20.4$   && $9557$ &\cr\tblrule
&& {\tt vvrk38} && $7.2E-9$ && $318$ && $1$    && $61.4$  && $37931$ &\cr\tblrule
&& {\tt vwdp5}  && $3.0E-9$ && $318$ && $718$  && $20.4$   && $11623$ &\cr\tblrule
&& {\tt vwrk38} && $4.6E-9$ && $318$ && $1835$ && $46.0$  && $34317$ &\cr\tblrule
\cr}}}
\par\bigskip

\end{exm}
\begin{exm}\label{Ex2}\rm
Here we take the coefficients matrix
$$A(t)=\alpha (\theta(t)-\sin(t)) \bmat{0 & 1 \cr -1 & 0} \,,$$
where $\theta(t) = \frac{\alpha}{1+\alpha^2}(\exp(-\alpha t)+\alpha\sin(t)-\cos(t))$
and we fixed $\alpha=100$.  Interval of integration is $[0,10]$
and exact solution is $Q(t)= \bmat{\cos(\theta(t))) & \sin(\theta(t)) 
\cr -\sin(\theta(t))  & \cos(\theta(t))}$. 
This is a difficult problem for all the methods considered, but all 
methods except
{\tt prk45} obtain accurate solutions as seen in Tables 4.3 and 4.4. 
The integrator {\tt RKF45} fails to integrate in one-step mode on 
the desired interval, though it succeeds on $[0,9]$, much less
efficiently than {\tt vtdp5, vvdp5, vwdp5}.
We expected the $\theta$ methods to run into difficulty in variable
stepsize mode, because of stiffness of the $\theta$ differential equation, 
but did not notice that.   

\bigskip
\centerline{
\vbox{\tabskip=0pt\offinterlineskip
\def\tblrule{\noalign{\hrule}}
\halign{ 
 \strut#&\vrule#\tabskip=1em plus2em
 &\hfil#\hfil&\vrule#    
 &\hfil#\hfil&\vrule#    
 &\hfil#\hfil&\vrule#    
 &\hfil#\hfil&\vrule#    
 \tabskip=0pt\cr         
\tblrule                      
&&\multispan{7}\hfil Table 4.3. Example \ref{Ex2}: fixed stepsize
$\Delta t=1.E-3$.
\hfil&\cr\tblrule
\tblrule\tblrule
&&\omit\hidewidth {\tt Meth} \hidewidth&&
  \omit\hidewidth {\tt err} \hidewidth&&
  \omit\hidewidth {\tt reimb} \hidewidth&&
  \omit\hidewidth {\tt cpu} \hidewidth&\cr\tblrule
\tblrule\tblrule
&& {\tt tdp5}  &&  $1.5E-12$  && $0$   && $169.$ &\cr\tblrule
&& {\tt trk38} &&  $1.5E-10$  && $0$   && $136.$ &\cr\tblrule
&& {\tt vdp5}  &&  $6.2E-12$  && $0$   && $209.$ &\cr\tblrule
&& {\tt vrk38} &&  $1.5E-10$  && $0$   && $170.$ &\cr\tblrule
&& {\tt wdp5}  &&  $6.2E-12$  && $0$   && $163.$ &\cr\tblrule
&& {\tt wrk38} &&  $1.6E-10$  && $0$   && $131.$ &\cr\tblrule
\cr}}}
\par\bigskip

\bigskip
\centerline{
\vbox{\tabskip=0pt\offinterlineskip
\def\tblrule{\noalign{\hrule}}
\halign{ 
 \strut#&\vrule#\tabskip=1em plus2em
 &\hfil#\hfil&\vrule#    
 &\hfil#\hfil&\vrule#    
 &\hfil#\hfil&\vrule#    
 &\hfil#\hfil&\vrule#    
 &\hfil#\hfil&\vrule#    
 &\hfil#\hfil&\vrule#    
 \tabskip=0pt\cr         
\tblrule                      
&&\multispan{11}\hfil Table 4.4. Example \ref{Ex2}: variable stepsize,
{\tt tol}$=1.E-8$.
\hfil&\cr\tblrule
\tblrule\tblrule
&&\omit\hidewidth {\tt Meth} \hidewidth&&
  \omit\hidewidth {\tt err} \hidewidth&&
  \omit\hidewidth {\tt reimb} \hidewidth&&
  \omit\hidewidth {\tt rejs} \hidewidth&&
  \omit\hidewidth {\tt cpu} \hidewidth&&
  \omit\hidewidth {\tt nsteps} \hidewidth&\cr\tblrule
\tblrule\tblrule
&& {\tt prk45}  && $-$ && $-$ && $-$    && $-$ && $-$ &\cr\tblrule
&& {\tt vtdp5}  && $5.3E-9$ && $0$ && $8$   && $1$   && $53$ &\cr\tblrule
&& {\tt vtrk38} && $5.1E-9$ && $0$ && $8$   && $3.0$ && $206$ &\cr\tblrule
&& {\tt vudp5}  && $1.5E-8$ && $0$ && $10$  && $2.0$ && $93$ &\cr\tblrule
&& {\tt vurk38} && $2.6E-8$ && $0$ && $1$   && $4.0$ && $280$ &\cr\tblrule
&& {\tt vvdp5}  && $6.9E-9$ && $0$ && $10$  && $2.9$ && $106$ &\cr\tblrule
&& {\tt vvrk38} && $5.9E-8$ && $0$ && $8$   && $5.0$ && $263$ &\cr\tblrule
&& {\tt vwdp5}  && $1.3E-8$ && $0$ && $12$  && $1.2$ && $66$ &\cr\tblrule
&& {\tt vwrk38} && $6.4E-9$ && $0$ && $20$  && $3.1$ && $238$ &\cr\tblrule
\cr}}}
\par\bigskip

\end{exm}
\begin{exm}\label{Ex3}\rm
This is a coefficients' matrix arising from a stiff two-point boundary
value problem having both boundary and interior layers.  We have
$$A(t) = \bmat{0 & 0 & 1 & 0 \cr t/(2\epsilon) & 0 & 1 & 1/2 \cr
1/\epsilon & 0 & 0 & 0 \cr
0 & 1/\epsilon & 1/\epsilon & -t/(2\epsilon)\cr}
$$
and we take $\epsilon=10^{-2}$. Interval of integration is $[-1,1]$.
Exact solution is not known.  The {\tt u}-methods fail to complete
the integration, but all other methods perform well.  For this problem,
see Tables 4.6-(i) and 4.6(ii), most rejections occur after having
computed all of $Q$.

\bigskip
\centerline{
\vbox{\tabskip=0pt\offinterlineskip
\def\tblrule{\noalign{\hrule}}
\halign{ 
 \strut#&\vrule#\tabskip=1em plus2em
 &\hfil#\hfil&\vrule#    
 &\hfil#\hfil&\vrule#    
 &\hfil#\hfil&\vrule#    
 \tabskip=0pt\cr         
\tblrule                      
&&\multispan{5}\hfil Table 4.5. Example \ref{Ex3}: fixed stepsize
$\Delta t=1.E-3$.
\hfil&\cr\tblrule
\tblrule\tblrule
&&\omit\hidewidth {\tt Meth} \hidewidth&&
  \omit\hidewidth {\tt reimb} \hidewidth&&
  \omit\hidewidth {\tt cpu} \hidewidth&\cr\tblrule
\tblrule\tblrule
&& {\tt tdp5}  && $2$ && $8.7$ &\cr\tblrule
&& {\tt trk38} && $2$ && $6.2$ &\cr\tblrule
&& {\tt vdp5}  && $3$ && $12.2$ &\cr\tblrule
&& {\tt vrk38} && $3$ && $8.8$ &\cr\tblrule
&& {\tt wdp5}  && $3$ && $10.0$ &\cr\tblrule
&& {\tt wrk38} && $3$ && $7.2$ &\cr\tblrule
\cr}}}
\par\bigskip

\bigskip
\centerline{
\vbox{\tabskip=0pt\offinterlineskip
\def\tblrule{\noalign{\hrule}}
\halign{ 
 \strut#&\vrule#\tabskip=1em plus2em
 &\hfil#\hfil&\vrule#    
 &\hfil#\hfil&\vrule#    
 &\hfil#\hfil&\vrule#    
 &\hfil#\hfil&\vrule#    
 &\hfil#\hfil&\vrule#    
 \tabskip=0pt\cr         
\tblrule                      
&&\multispan{9}\hfil Table 4.6-(i). Example \ref{Ex3}: variable stepsize,
{\tt tol}$=1.E-8$.
\hfil&\cr\tblrule
\tblrule\tblrule
&&\omit\hidewidth {\tt Meth} \hidewidth&&
  \omit\hidewidth {\tt reimb} \hidewidth&&
  \omit\hidewidth {\tt rejs/first} \hidewidth&&
  \omit\hidewidth {\tt cpu} \hidewidth&&
  \omit\hidewidth {\tt nsteps} \hidewidth&\cr\tblrule
\tblrule\tblrule
&& {\tt prk45}  && $-$   && $-$       && $1.4$  && $252$ &\cr\tblrule
&& {\tt vtdp5}  && $2$   && $11/2$    && $1$    && $221$ &\cr\tblrule
&& {\tt vtrk38} && $2$   && $15/4$    && $1.7$  && $628$ &\cr\tblrule
&& {\tt vudp5}  && $-$   && $-$       && $-$    && $-$ &\cr\tblrule
&& {\tt vurk38} && $-$   && $-$       && $-$    && $-$ &\cr\tblrule
&& {\tt vvdp5}  && $3$   && $9/1$     && $1.3$  && $217$ &\cr\tblrule
&& {\tt vvrk38} && $3$   && $14/2$    && $2.7$  && $612$ &\cr\tblrule
&& {\tt vwdp5}  && $3$   && $10/2$    && $1.1$  && $228$ &\cr\tblrule
&& {\tt vwrk38} && $3$   && $13/3$    && $2.3$  && $649$ &\cr\tblrule
\cr}}}
\par\bigskip

\bigskip
\centerline{
\vbox{\tabskip=0pt\offinterlineskip
\def\tblrule{\noalign{\hrule}}
\halign{ 
 \strut#&\vrule#\tabskip=1em plus2em
 &\hfil#\hfil&\vrule#    
 &\hfil#\hfil&\vrule#    
 &\hfil#\hfil&\vrule#    
 &\hfil#\hfil&\vrule#    
 \tabskip=0pt\cr         
\tblrule                      
&&\multispan{7}\hfil Table 4.6-(ii). Example \ref{Ex3}: column-rejections,
{\tt tol}$=1.E-8$.
\hfil&\cr\tblrule
\tblrule\tblrule
&&\omit\hidewidth {\tt Meth} \hidewidth&&
  \omit\hidewidth {\tt 1st} \hidewidth&&
  \omit\hidewidth {\tt 2nd} \hidewidth&&
  \omit\hidewidth {\tt 3rd} \hidewidth&\cr\tblrule
\tblrule\tblrule
&& {\tt vtdp5}  &&  2 &&  1 && 8  &\cr\tblrule
&& {\tt vtrk38} &&  4 &&  0 && 11 &\cr\tblrule
&& {\tt vvdp5}  &&  1 &&  0 && 8 &\cr\tblrule
&& {\tt vvrk38} &&  2 &&  0 && 12 &\cr\tblrule
&& {\tt vwdp5}  &&  2 &&  1 && 7 &\cr\tblrule
&& {\tt vwrk38} &&  3 &&  0 && 10 &\cr\tblrule
\cr}}}
\par\bigskip

\end{exm}
\begin{exm}\label{Ex4}\rm  This is a fairly simple problem quite useful 
for testing purposes.  We have the coefficient matrix
$$A(t)=Q(t)D(t)Q^T(t)+\dot Q(t)Q^T(t)\,,$$
where $D(t)=\diag(1,\cos (t),-\frac{1}{2\sqrt{t+1}}, -10)$, 
$Q(t)=\bmat{1 & 0 & 0 \cr 0 & Q_\beta(t) & 0 \cr 0 & 0 & 1}
\bmat{Q_\alpha(t) & 0 \cr 0 & Q_\alpha(t)}$, and
$Q_\gamma(t)=\bmat{\cos (\gamma t) & \sin (\gamma t) \cr
-\sin (\gamma t) & \cos (\gamma t)}$.  We report on results for
$\alpha=1$ and $\beta=\sqrt{2}$ and integrate on $[0,100]$.  
Notice the inadequacy of the $u$-variables. 
Such inadequacy is hidden in
constant stepsize mode, but it becomes apparent in variable stepsize
mode: the $u$-vector becomes very poorly scaled.
In variable stepsize, nearly all 
rejections occur while triangularizing the first column.
From Table 4.8, observe how {\tt prk45} takes 
more steps and is less accurate than {\tt vtdp5, vvdp5, vwdp5}.

\bigskip

\centerline{
\vbox{\tabskip=0pt\offinterlineskip
\def\tblrule{\noalign{\hrule}}
\halign{ 
 \strut#&\vrule#\tabskip=1em plus2em
 &\hfil#\hfil&\vrule#    
 &\hfil#\hfil&\vrule#    
 &\hfil#\hfil&\vrule#    
 &\hfil#\hfil&\vrule#    
 \tabskip=0pt\cr         
\tblrule                      
&&\multispan{7}\hfil Table 4.7. Example \ref{Ex4}: fixed stepsize
$\Delta t=1.E-3$.
\hfil&\cr\tblrule
\tblrule\tblrule
&&\omit\hidewidth {\tt Meth} \hidewidth&&
  \omit\hidewidth {\tt err} \hidewidth&&
  \omit\hidewidth {\tt reimb} \hidewidth&&
  \omit\hidewidth {\tt cpu} \hidewidth&\cr\tblrule
\tblrule\tblrule
&& {\tt tdp5}  && $1.6E-10$ && $27$  && $22.7$ &\cr\tblrule
&& {\tt trk38} && $1.5E-10$ && $27$  && $16.9$  &\cr\tblrule
&& {\tt vdp5}  && $1.6E-10$ && $77$  && $24.3$ &\cr\tblrule
&& {\tt vrk38} && $1.5E-10$ && $77$  && $19.1$ &\cr\tblrule
&& {\tt wdp5}  && $1.6E-10$ && $77$  && $22.2$ &\cr\tblrule
&& {\tt wrk38} && $1.5E-10$ && $77$  && $17.3$ &\cr\tblrule
\cr}}}
\par\bigskip

\bigskip
\centerline{
\vbox{\tabskip=0pt\offinterlineskip
\def\tblrule{\noalign{\hrule}}
\halign{ 
 \strut#&\vrule#\tabskip=1em plus2em
 &\hfil#\hfil&\vrule#    
 &\hfil#\hfil&\vrule#    
 &\hfil#\hfil&\vrule#    
 &\hfil#\hfil&\vrule#    
 &\hfil#\hfil&\vrule#    
 &\hfil#\hfil&\vrule#    
 \tabskip=0pt\cr         
\tblrule                      
&&\multispan{11}\hfil Table 4.8. Example \ref{Ex4}: variable stepsize,
{\tt tol}$=1.E-8$.
\hfil&\cr\tblrule
\tblrule\tblrule
&&\omit\hidewidth {\tt Meth} \hidewidth&&
  \omit\hidewidth {\tt err} \hidewidth&&
  \omit\hidewidth {\tt reimb} \hidewidth&&
  \omit\hidewidth {\tt rejs/first} \hidewidth&&
  \omit\hidewidth {\tt cpu} \hidewidth&&
  \omit\hidewidth {\tt nsteps} \hidewidth&\cr\tblrule
\tblrule\tblrule
&& {\tt prk45}  && $2.1E-7$  && $-$   && $-$  && $1.1$ && $5053$ &\cr\tblrule
&& {\tt vtdp5}  && $7.7E-9$  && $27$   && $94/89$    && $1$ && $4533$ &\cr\tblrule
&& {\tt vtrk38} && $1.2E-8$  && $27$   && $107/100$  && $2.1$   && $13010$ &\cr\tblrule
&& {\tt vudp5}  && $-$       && $-$   && $-$         && $-$ && $-$ &\cr\tblrule
&& {\tt vurk38} && $-$       && $-$   && $-$         && $-$ && $-$ &\cr\tblrule
&& {\tt vvdp5}  && $1.2E-8$  && $77$  && $80/80$     && $1.1$   && $3967$ &\cr\tblrule
&& {\tt vvrk38} && $1.2E-8$  && $77$  && $99/99$     && $2.1$  && $11169$ &\cr\tblrule
&& {\tt vwdp5}  && $1.4E-8$  && $77$  && $103/92$    && $1.1$   && $4370$ &\cr\tblrule
&& {\tt vwrk38} && $2.8E-8$  && $77$  && $125/115$   && $2.3$  && $12694$ &\cr\tblrule
\cr}}}
\par\bigskip
\end{exm}

\begin{exm}\label{Ex7}\rm
Here we take the diagonal coefficient matrix
$$A(t)=\diag(-\frac{1}{2\sqrt{t+1}}, -10,\cos(t),1)\,,$$
which is just a permutation of $D(t)$ in Example \ref{Ex4}.
We take both $Q(0)=I$ and a random initial condition $Q(0)$, to verify
if the diagonal elements of $D$ will become reordered.  
When $Q(0)=I$, then the computed solution remains $Q(t)=I$ for all $t$.
For randomly chosen $Q(0)$, the change
of variables $Q$ does instead reorder the diagonal of $\tilde A$:
$\tilde A_{11}\geq \cdots \geq \tilde A_{44}$ (see Figure 4.1;
results obtained with {\tt vtdp5}).

\begin{figure}
\vskip-2.0in
\centerline{
\epsfysize=7.0in
\hskip-.3in\epsfbox{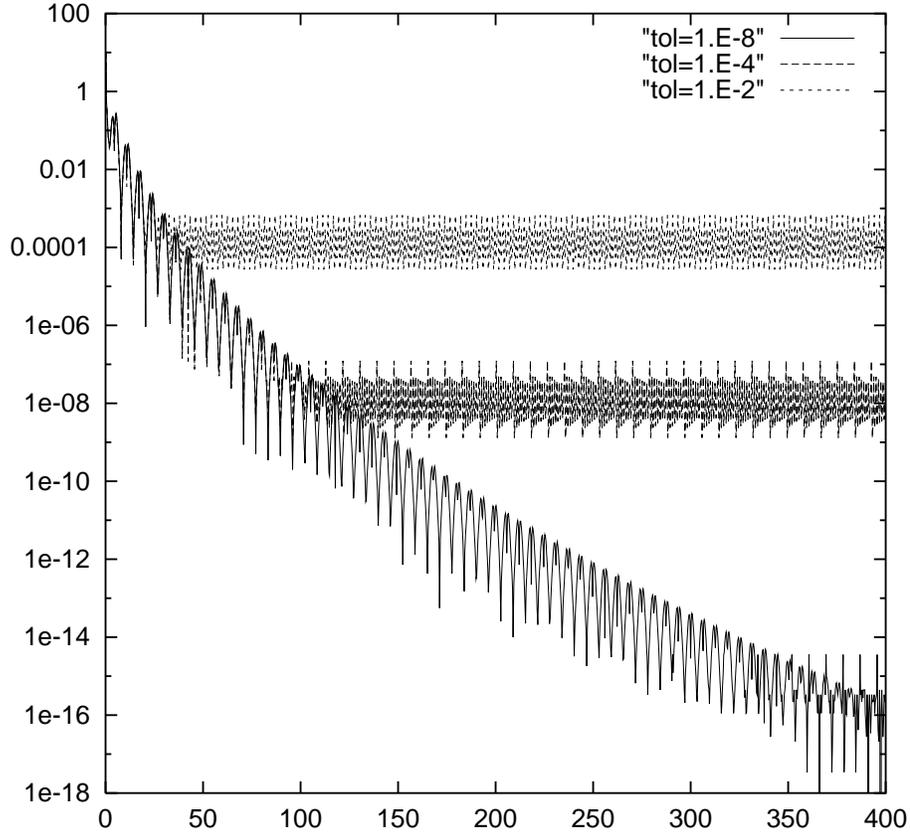}}
\caption{Plot of $t$ vs. $\log_{10}(||D(t)-\diag(\tilde A(t))||_\infty)$
for {\tt vtdp5} applied to Example \ref{Ex7} for different tolerances using
a ``random'' $Q(0)$.} 
\end{figure}

\bigskip

\end{exm}

\begin{exm}\label{Ex6}\rm
This is an example of a constant coefficients matrix.
The equation (\ref{QRflow}) reduces to a so-called isospectral flow.  
For these problems, it is known that the diagonal of the
upper $(p,p)$ block of the transformed matrix must
eventually converge to the real parts of the leading
eigenvalues of $A$.
Far from advocating the ideas set forth in this paper as a mean to solve
the eigenvalue problem, we included this problem 
because we wanted to highlight 
(on a problem of arbitrarily high dimensions) how most stepsize
rejections occur ahead of having found all of $Q$; for this reason, 
we report on results obtained with the variable stepsize codes 
{\tt vtdp5, vvdp5, vwdp5}, and, for comparison, with {\tt prkf45},
in spite of the fact that probably constant coefficients
problems can be efficiently solved also with constant stepsize.
Also, we wanted to see how ill conditioning (i.e., poor separation)
of the eigenvalues affected convergence.  

The coefficient matrix is the upper Hessenberg Frank matrix:
$$A=\pmatrix{n & n-1 & n-2 & \dots & 1 \cr
n-1 & n-1 & n-2 & \dots & 1 \cr
0 & n-2 & n-2 & \dots & 1 \cr
\vdots & \ddots & \ddots & \ddots & \vdots \cr
0 & \dots & 0 & 1 & 1} \,.$$
We take $n=25$, and seek $Q$ triangularizing 
$X\in \R^{n\times p}$, $X'=AX, \ X(0)=\bmat{I_p \cr 0}$, with
$p=13$.  Interval of integration is $[0,b]$.  Exact solution is 
not known, so we report on the defect of
the diagonal of the upper $(p,p)$ block of the 
transformed matrix $\tilde A(b)$  from the real parts of the 
leading eigenvalues; this we list as {\tt errd}.  
At four digits, the eigenvalues are 
$77.9837$, $60.5984$, $47.7777$, $37.5667$,
$29.2021$, $22.2856$, $16.5772$, $11.9193$,
$8.2006$, $5.3359$, $3.2479$,  
$1.8495$, $1.1841 \pm 0.2634i$, $0.8147 \pm 0.6124i$,
$0.4098 \pm 0.7755i$, $0.0051 \pm 0.7737i$, $-0.3373 \pm 0.6031i$,
$-0.5452 \pm 0.3177i$, $-0.6070$.

In the next table, we report on selected runs with the variable
stepsize codes for $b=100$.  The value of {\tt errd} is entirely
due to lack of convergence to the real part of the 13-th eigenvalue, 
$1.1841$, all larger eigenvalues
having already converged to several digits.  Notice that (for most cases)
{\tt errd} does not change appreciably by decreasing the tolerance value, 
which indicates that only by increasing the length of integration one would 
be able to approximate the 13th eigenvalue more accurately.  Also, we notice
that our codes are clearly superior to {\tt prkf45}.  Finally, in our codes,
a large majority of rejections occur while integrating for the 
4th or 5th columns; for example, for {\tt vtdp5} with {\tt tol=1.E-6}, the
rejection patterns is as follows: $[1,0,11,372,104,4,0,0,1,14,1,4,3]$, where
each number in this vector refers to rejections occurred while
triangularizing columns $1$ to $13$.

\bigskip
\centerline{
\vbox{\tabskip=0pt\offinterlineskip
\def\tblrule{\noalign{\hrule}}
\halign{ 
 \strut#&\vrule#\tabskip=1em plus2em
 &\hfil#\hfil&\vrule#    
 &\hfil#\hfil&\vrule#    
 &\hfil#\hfil&\vrule#    
 &\hfil#\hfil&\vrule#    
 &\hfil#\hfil&\vrule#    
 &\hfil#\hfil&\vrule#    
 &\hfil#\hfil&\vrule#    
 \tabskip=0pt\cr         
\tblrule                      
&&\multispan{13}\hfil Table 4.9. Example \ref{Ex6}: variable stepsize.
\hfil&\cr\tblrule
\tblrule\tblrule
&&\omit\hidewidth {\tt Meth} \hidewidth&&
  \omit\hidewidth {\tt tol} \hidewidth&&
  \omit\hidewidth {\tt errd} \hidewidth&&
  \omit\hidewidth {\tt reimb} \hidewidth&&
  \omit\hidewidth {\tt rejs} \hidewidth&&
  \omit\hidewidth {\tt cpu} \hidewidth&&
  \omit\hidewidth {\tt nsteps} \hidewidth&\cr\tblrule
\tblrule\tblrule
&& {\tt prk45}  && $1.E-4$ && $1.8E-1$ && -    && -         && $1.7$  && $5365$ &\cr\tblrule
&& {\tt vtdp5}  && $1.E-4$ && $9.6E-2$ && $73$ && $525$ && $1.5$  && $2391$ &\cr\tblrule
&& {\tt vvdp5}  && $1.E-4$ && $1.1E0$ && $77$ && $516$ && $1.2$  && $2459$ &\cr\tblrule
&& {\tt vwdp5}  && $1.E-4$ && $3.0E0$ && $85$ && $504$ && $1$  && $2462$ &\cr\tblrule
&& {\tt prk45}  && $1.E-6$ && $1.8E-1$ && -    && -         && $1.7$  && $5430$ &\cr\tblrule
&& {\tt vtdp5}  && $1.E-6$ && $1.8E-1$ && $65$ && $515$ && $1.6$  && $2459$ &\cr\tblrule
&& {\tt vvdp5}  && $1.E-6$ && $5.6E-2$ && $35$ && $501$ && $1.2$  && $2491$ &\cr\tblrule
&& {\tt vwdp5}  && $1.E-6$ && $1.1E-1$ && $34$ && $476$ && $1.$  && $2481$ &\cr\tblrule
\cr}}}
\par\bigskip

\end{exm}

\section{Conclusions}\label{s5}

The purpose of this work has been manifold.  We discussed general
design choices people should have in mind when devising techniques
for solving (\ref{QRflow}).  We reinterpreted methods based
on Householder and Givens transformations as methods which give
a trajectory on the smooth manifold of orthonormal matrices, and
use overlapping local charts to parameterize the manifold.  We
gave new formulations and new implementations of methods based on 
Householder and Givens transformations.  Finally, we presented a suite
of {\tt FORTRAN} codes for solving (\ref{QRflow}) by our techniques.
We believe that the methods put forward in this work, and their
implementation as laid down here, are a very sensible way to solve
(\ref{QRflow}), and should not be ignored by anyone interested
in comparative performance of other techniques for solving (\ref{QRflow}).
In this spirit, we welcome other people to use our 
codes\footnote{e-mail to either author}.

Our study showed that all design scopes we had set at the beginning
are satisfied by appropriate implementation of our techniques.
There is certainly room for improving our codes, and we anticipate
some work towards more efficient implementation, in particular insofar
as data structure, memory requirement, and extensive use of the {\tt BLAS}
libraries.  However, we believe that our implementations
are now sophisticated enough that some conclusions and recommendations
for future work can be given.
\begin{enumerate}
\item Much as one should have expected, in variable stepsize the
{\tt dp5} codes outperform the {\tt rk38} codes.  In constant
stepsize, however, the {\tt rk38} codes are less expensive and
this makes up for their reduced accuracy.
\item  A striking difference with the standard linear algebra situation
occurs when using Householder methods.  There, different representations
of the Householder vectors are chiefly a matter of storage efficiency.
In our time dependent setting, however, 
Householder methods in the $u$-variables should be discarded:
they are prone to instabilities,
and this is clearly betrayed in a variable stepsize implementation.
Instead, in the
$v$-variables, integrated with a local projection step, and - especially -
the $w$-variables, Householder methods are robust and also handle 
well large problems.  
\item Givens methods based on
the $\theta$-variables revealed very accurate and efficient.
However, the frequent need to compute sines and cosines and inverse
trigonometric functions is a potential
drawback with this approach.  A careful implementation based on
direct integration of the $(\cos,\sin)$ values has not been carried out, 
but could be an interesting endeavor.
\item The implementation in variable stepsize has proved very
rewarding, especially since it confirmed our expectation 
that often most rejections occur far ahead of having computed all of $Q$.  
\item  It would be interesting to implement our methods by
exploiting the rewriting (\ref{transQ}), and hence
solving (\ref{QRflow2}) always starting near $\bmat{I_p \cr 0}$.  
This may prove beneficial for all methods, also
those based on the $u$-variables.
\item  With the present level of implementation, and
all things considered, the methods based on $\theta$ and $w$ variables
are probably the best, followed by the $v$-methods.
\item The relative comparison with the projected integrator {\tt prk45} 
appears favorable to our new codes.  In particular, our codes generally
require fewer steps, and are more accurate and less expensive.  
Moreover, the {\tt vtdp5, vvdp5, vwdp5} 
codes never failed to complete the integration; instead,
{\tt prk45}, in spite of the certainly more sophisticated implementation
of the integrator {\tt RKF45}, was occasionally unable to complete 
the integration.
\end{enumerate}


\end{document}